\documentclass[12pt]{amsart}
\usepackage{amsmath,amssymb,latexsym,psfig,colordvi}

\newcommand{\sgn}{{\rm sgn}}

\newcommand{\fs}{\mathfrak S}
\newcommand{\fg}{\mathcal G}



\newtheorem{dfn}[equation]{Definition}
\newtheorem{alg}[equation]{Algorithm}
\newtheorem{thm}[equation]{Theorem}
\newtheorem{remark}[equation]{Remark}
\newtheorem{cor}[equation]{Corollary}
\newtheorem{lem}[equation]{Lemma}
\newtheorem{example}[equation]{Example}
\newenvironment{rem}{\begin{remark}\rm}{\end{remark}}
\newenvironment{ex}{\begin{example}\rm}{\end{example}}

\makeatletter\@addtoreset{equation}{section}\makeatother

\newcommand{\bC}{{\mathbb C}}
\newcommand{\spa}{\;\: }
 
\newcommand{\bZ}{{\mathbb Z}}

\newcommand{\eg}{{\rm end}(\gamma)}

\newcommand{\spandsp}{\mbox{$\spa\;\text{and}\;\spa$}}
\newcommand{\s}{\sigma}

\newcounter{FNC}[page]
\def\fauxfootnote#1{{\addtocounter{FNC}{2}$^\fnsymbol{FNC}$%
     \let\thefootnote\relax\footnotetext{$^\fnsymbol{FNC}$#1}}}

\title{A Pieri-type formula for the ${K}$-theory of a flag manifold}
\author{Cristian Lenart}
\address{Department of Mathematics and Statistics\\
         State University of New York at Albany\\
         Albany, NY \ 12222\\
         USA}
\email{lenart@csc.albany.edu}
\urladdr{http://math.albany.edu/math/pers/lenart}

\author{Frank Sottile}
\address{Department of Mathematics\\
         Texas A\&M University\\
         College Station\\
         TX \ 77843\\
         USA}
\email{sottile@math.tamu.edu}
\urladdr{http://www.math.tamu.edu/\~{}sottile}


\thanks{Research of Lenart supported by SUNY Albany Faculty Research Award 1039703}
\thanks{Research of Sottile supported in part by the Clay Mathematical Institute,
  the MSRI, and NSF CAREER grant DMS-0134860}

\begin{document}

\begin{abstract}
 We derive explicit Pieri-type multiplication formulas in the Grothendieck ring
 of a flag variety.
 These expand the product of an arbitrary Schubert class and a special Schubert
 class in the basis of Schubert classes. 
 These special Schubert classes are indexed by a cycle which has either the form 
 $(k{-}p{+}1,k{-}p{+}2,\ldots,k{+}1)$ or the form $(k{+}p,k{+}p{-}1,\ldots,k)$,
 and are pulled back from a Grassmannian projection. 
 Our formulas are in terms of certain labeled chains in the $k$-Bruhat order on
 the symmetric group and are combinatorial in that they involve no cancellations. 
 We also show that the multiplicities in the Pieri formula are naturally certain
 binomial coefficients.
\end{abstract}

\maketitle

\vspace{-5mm}

\section*{Introduction} 

Classically, Schubert calculus is concerned with enumerative problems in
geometry, such as counting the lines satisfying some generic
intersection conditions. 
This enumeration is accomplished via a calculation in the
cohomology ring of the space of potential solutions such as a
Grassmannian~\cite{KL72}. 
The cohomology ring of a Grassmannian is well-understood combinatorially
through the Littlewood-Richardson rule.
Less understood, particularly in combinatorial terms, are extensions to more
general flag varieties and to more general cohomology theories, such as
equivariant cohomology, quantum cohomology, or $K$-theory.
The ``modern Schubert calculus'' is concerned with the geometry and
combinatorics of these extensions. 
We study the multiplicative structure of the
Grothendieck ring ($K$-theory) of the manifold of flags in $n$-space, giving a
Pieri-type formula in the sense of~\cite{sotprf}.
Our formulas and their proofs are based on combinatorics of the Bruhat order
on the symmetric group, and they highlight new properties of this order. 

The flag variety has an algebraic Schubert cell
decomposition.
Consequently, classes of structure sheaves of Schubert varieties (Schubert
classes) form an integral basis of its Grothendieck ring, which is 
indexed by permutations. 
A major open problem in the modern Schubert calculus is to determine the
$K$-theory Schubert structure constants, which express a product of two
Schubert classes in terms of this Schubert basis.
Brion~\cite{brion} proved that these coefficients alternate in sign in a
specified manner.

In the passage from the filtered Grothendieck ring to its associated graded
ring, which is isomorphic to the cohomology ring, 
Schubert classes are mapped to classes of Schubert varieties (also called
Schubert classes).
In this way, the cohomology Schubert structure constants are certain 
$K$-theory Schubert structure constants. 
These cohomology constants are Littlewood-Richardson
constants when the Schubert classes come from a Grassmannian.
When one class is of a hypersurface Schubert variety,
Monk~\cite{mongfm} gave a formula for these constants.
Monk's formula highlights the importance of a suborder of the Bruhat order on
the symmetric group called the $k$-Bruhat order.
For example, the Pieri formula in~\cite{sotprf} uses chains in 
the $k$-Bruhat order to express the multiplication of a Schubert class by a
special Schubert class pulled back from the Grassmannian of
$k$-planes in $\bC^n$.

Buch~\cite{buch} gave the $K$-theory formula for the
product of two classes pulled back from the same Grassmannian.
Until now, the only general formula in the Grothendieck ring of the flag
variety is the analog of Monk's formula for multiplication by the structure
sheaf of a hypersurface Schubert variety~\cite{falpfg,GR04,lenktv,LP03,parpcf}.
The formula in~\cite{lenktv} is in terms of chains in the $k$-Bruhat order.
We give a Pieri-type formula in the Grothendieck ring, which 
generalizes both the $K$-theory Monk formula of~\cite{lenktv} and the
cohomology Pieri formula~\cite{sotprf}.
The formula is in terms of chains in the $k$-Bruhat order, but with some
covers marked.
The unmarked covers satisfy a condition from the $K$-theory Monk formula, while
the marked covers satisfy a condition from the cohomology Pieri formula.

We work in the algebraic-combinatorial theory of Grothendieck polynomials.
These distinguished polynomial representatives of $K$-theory Schubert classes
were introduced by Lascoux and Sch\"utzenberger~\cite{lassfm} and studied
further in~\cite{lasagv}.
The transition formula~\cite{lastgp,lenktv} gives a recursive construction of
Grothendieck polynomials, and the recursion for polynomials representing
special Schubert classes is the basis of our proof.
The other main ingredient of our proof is a Monk-like formula for multiplying a
Grothendieck polynomial by a variable given in~\cite{lenktv}.
This recursion was suggested as a basis for a proof of the Pieri
formula in cohomology by Lascoux and Sch\"utzenberger~\cite{lasps}.
Manivel presented such a proof in his book~\cite[p.~94]{manfsp}, but that
proof contains a subtle error, omitting some important subcases. 
We correct that omission in our proof, see Remark~\ref{R:M3}.\smallskip

In Section~1, we give basic definitions and background concerning Grothendieck
polynomials, describe the Monk formulas in $K$-theory and the Pieri formula
in cohomology that this work generalizes, and state our Pieri-type formula.
A uniqueness result about chains proven in Section 2 implies a version of the
formula in which the coefficients are naturally certain binomial coefficients.
Section 2 also collects some technical results on the
Bruhat order used in our proof of the Pieri-type formula, which
occupies Section 3.
We conclude in Section 4 with a dual Pieri-type formula and a discussion of the
specialization of our formula to the Pieri-type formula in the $K$-theory of a
Grassmannian~\cite{lencak}. 

\section{Grothendieck polynomials and the Pieri formula}

We first introduce Schubert and Grothendieck polynomials.
For more information, see~\cite{fulyt,lasagv,macnsp,manfsp}. 
We next state the known Monk and Pieri-type formulas in the cohomology and the 
$K$-theory of a flag manifold,
and then state our Pieri-type formula.

\subsection{Schubert and Grothendieck polynomials}
Let $Fl_n$ be the variety of complete flags 
$(\{0\}=V_0\subset V_1\subset\ldots\subset V_n=\bC^n)$ in $\bC^n$. 
This algebraic manifold has dimension $\binom{n}{2}$.  
Its integral cohomology ring
$H^*(Fl_n)$ is isomorphic to $\bZ[x_1,\ldots,x_n]/I_n$, where the ideal $I_n$ is
generated by the nonconstant symmetric polynomials in
$x_1,\ldots,x_n$, and $x_i$ has cohomological degree 2. 
For this, the element $x_i$ is identified with the Chern class of the 
dual $L^*_i$ to the tautological line bundle $L_i:=V_i/V_{i-1}$.
The variety $Fl_n$ is a disjoint union of cells indexed by permutations $w$ in
the symmetric group $S_n$. 
The closure of the cell indexed by $w$ is the {\em Schubert variety} $X_w$, which
has codimension $\ell(w)$, the length of $w$ or the number of its inversions. 
The Schubert polynomial $\fs_w(x)$ (defined below) is
a certain polynomial representative for the cohomology class corresponding to
$X_w$. It is a homogeneous polynomial in $x_1,\ldots,x_{n-1}$ of degree $\ell(w)$
with nonnegative integer coefficients.  

The Grothendieck group $K^0(Fl_n)$ of complex vector bundles on $Fl_n$ 
is isomorphic to its Grothendieck group of coherent sheaves.
As abstract rings,  $K^0(Fl_n)$ and $H^*(Fl_n)$ are isomorphic.
Here, the variable $x_i$ is the $K$-theory Chern class $1{-}1/y_i$ of the line
bundle $L_i^*$, where $y_i$ represents $L_i$ in the Grothendieck
ring. 
The classes of the structure sheaves of Schubert varieties form a 
basis of $K^0(Fl_n)$. 
The class indexed by $w$ is represented by the Grothendieck polynomial
$\fg_w(x)$. 
This inhomogeneous polynomial in $x_1,\ldots,x_{n-1}$ has lowest degree 
homogeneous component equal to the Schubert polynomial $\fs_w(x)$. 

The construction of Schubert and Grothendieck polynomials is based on the {\em
  divided difference operators} $\partial_i$ and the {\em isobaric divided
  difference operators} $\pi_i$.
As operators on $\bZ[x_1,x_2,\ldots]$, these are defined as follows: 
 \begin{equation}\label{defpi}
  \partial_i\ :=\ \frac{1-s_i}{x_i-x_{i+1}}\quad\mbox{and}\quad
  \pi_i\ :=\ \partial_i(1-x_{i+1})=1+(1-x_i)\partial_i\,.
\end{equation}
Here $s_i$ is the transposition $({i,i+1})$, which interchanges the
variables $x_i$ and $x_{i+1}$, $1$ is the identity operator, and $x_i$ is 
multiplication by the corresponding variable.  

Schubert and Grothendieck polynomials are defined inductively for each
permutation $w$ in $S_n$ by setting
$\fs_{\omega_0}(x)=\fg_{\omega_0}(x)=x_1^{n-1}x_2^{n-2}\dotsb
x_{n-1}$ where $\omega_0:=n\dotsc 21$ is the longest permutation in $S_n$
(in one-line notation), and by letting
 \begin{equation}\label{rec1}
   \fs_{ws_i}(x)\ :=\ \partial_i\,\fs_w(x)\,\quad\mbox{and}\quad
   \fg_{ws_i}(x)\ :=\ \pi_i\,\fg_w(x)\,,\quad\mbox{if}\quad 
   \ell(ws_i)\ =\ \ell(w)-1\,. 
 \end{equation}
A Grothendieck polynomial does not depend on the chosen chain
in the weak order on $S_n$ from $\omega_0$ to $w$ because the operators $\pi_i$
satisfy the braid relations 
 \begin{align*}
  \pi_i\pi_j&=\pi_j\pi_i\qquad\qquad\mbox{if}\spa |i-j|\ge 2\,,\\
  \pi_i\pi_{i+1}\pi_i&=\pi_{i+1}\pi_i\pi_{i+1}\,,
 \end{align*}
and similarly for Schubert polynomials. 
While defined for $w\in S_n$, the Schubert and Grothendieck polynomials
$\fs_w(x)$ and $\fg_w(x)$ do not not depend on $n$.
Thus we may define them for $w$ in $S_{\infty}$, where 
$S_{\infty}:=\bigcup_n S_n$ under the usual inclusion 
$S_n\hookrightarrow S_{n+1}$. 
Both the Schubert polynomials $\fs_w(x)$ and the
Grothendieck polynomials $\fg_w(x)$ form bases of
$\bZ[x_1,x_2,\ldots]$, as $w$ ranges over $S_{\infty}$.

\subsection{Known Monk and Pieri-type formulas}
The covering relations in the Bruhat order are $v\lessdot w=v({a,b})$, where
$\ell(w)=\ell(v)+1$. 
We denote this by 
 \[
    v\xrightarrow{\,({a,b})\,}w\,.
 \]
A permutation $v$ admits a cover $v\lessdot v(a,b)$ with $a<b$ and
$v(a)<v(b)$ if and only if whenever $a<c<b$, then either $v(c)<v(a)$ or else
$v(b)<v(c)$. Call this the {\it cover condition}. 
The $k$-Bruhat order first appeared in the context of Monk's
formula~\cite{mongfm}, and was studied in more detail in~\cite{basspb,BS99a}.
It is the suborder of the Bruhat order where the covers are restricted to 
those $v\lessdot v(a,b)$ with $a\le k<b$. 
Throughout this paper, we will consider the Bruhat order and the $k$-Bruhat order on $S_\infty$ rather than on a particular symmetric group $S_n$. We will use the following order on pairs of positive integers to compare covers
in a $k$-Bruhat order:
 \begin{equation}\label{ord}
  (a,b)\prec (c,d)\quad\mbox{if and only if}\quad  (b>d)\ \mbox{or}\ 
  (b=d\ \mbox{and}\ a< c)\,. 
 \end{equation}
The Monk formula for Grothendieck polynomials is a formula for multiplication by 
$\fg_{s_k}(x)$.

\begin{thm}\label{kmonk}\cite{lenktv} 
  We have that
 \[
    \fg_v(x)\,\fg_{s_k}(x)\ =\ \sum_{\gamma}(-1)^{\ell(\gamma)-1}\fg_{\eg}(x)\,,
 \]
 where the sum is over all saturated chains $\gamma$ in the $k$-Bruhat order (on $S_\infty$) 
 \[
   v=v_0\xrightarrow{\,({a_1,b_1})\,} 
                     v_1\xrightarrow{\,({a_2,b_2})\,}\  
       \dotsb\ 
     \xrightarrow{\,({a_p,b_p})\,}v_p\ =\ \eg\,,
 \]
 with $p=\ell(\gamma)\geq 1$, and 
 \begin{equation}\label{E:inc_chain}
    (a_1,b_1)\prec(a_2,b_2)\prec\dotsb \prec (a_p,b_p)\,.
 \end{equation}
 This formula has no cancellations and is multiplicity free, which
 means that the coefficients in the right-hand side are $\pm 1$ after
 collecting terms.   
\end{thm}

The proof of Theorem~\ref{kmonk} is based on the following formula for
multiplying an arbitrary Grothendieck polynomial by a single
variable.

\begin{thm}\label{xkmonk}\cite{lenktv} 
 We have that
  \begin{equation}\label{xk}
      x_k\,\fg_v(x)\ =\ \sum_{\gamma} \s(\gamma)\fg_{\eg}(x)\,,
  \end{equation}
 where the sum is over all saturated chains $\gamma$ in the Bruhat order (on $S_\infty$) 
 \[
   v=v_0\xrightarrow{\,({a_1,k})\,} 
                    v_1\xrightarrow{\,({a_2,k})\,}\ 
           \dotsb\ 
         \xrightarrow{\,({a_p,k})\,}v_p 
         \xrightarrow{\,({k,b_1})\,}v_{p+1}
         \xrightarrow{\,({k,b_2})\,}\ \dotsb\ 
         \xrightarrow{\,({k,b_q})\,}v_{p+q} 
       = \eg\,,
  \]
 where $p,q\ge 0$, $p+q\ge 1$, $\s(\gamma):=(-1)^{q+1}$, and 
  \begin{equation}\label{Eq:Monk_Condition}
     a_p<a_{p-1}<\dotsb<a_1<k<b_q<b_{q-1}<\dotsb<b_1\,.
  \end{equation}
\end{thm}

The Pieri formula for Schubert polynomials expresses the product of a
Schubert polynomial with an elementary symmetric polynomial
$e_p(x_1,\ldots,x_k)$, which is the Schubert polynomial indexed by the cycle
$c_{[k,p]}:=(k{-}p{+}1,k{-}p{+}2,\ldots,k{+}1)$, where $k\ge p\ge 1$.
There is a similar formula for multiplication by the 
homogeneous symmetric polynomial $h_p(x_1,\ldots,x_k)$, which is the
Schubert polynomial indexed by the cycle $r_{[k,p]}:=(k+p,k+p-1,\ldots,k)$.  
More generally, the Schur polynomial $s_\lambda(x_1,\ldots,x_k)$ indexed by the
partition $\lambda=(\lambda_1\ge \lambda_2\ge\dotsb\ge\lambda_k\geq0)$ coincides
with the Schubert polynomial indexed the permutation $v$ with a unique descent
at $k$, where $v(i)=\lambda_{k+1-i}+i$ for $1\le i\le k$.
We denote this by $v(\lambda,k)$ and call it a {\em Grassmannian permutation}.
The corresponding Schubert classes are pulled back from the projection
of $Fl_n$ to the Grassmannian of $k$-planes.

\begin{thm}\label{schub}\cite{lasps,manfsp, posqvp, sotprf} 
 We have that
 \[
    \fs_v(x)\,e_{p}(x_1,\ldots,x_k)=\sum_\gamma\fs_{\eg}(x)\,,
 \]
 where the sum is over all saturated chains $\gamma$ in the $k$-Bruhat order (on $S_\infty$) 
 \[
   v=v_0\xrightarrow{\,({a_1,b_1})\,}  
                    v_1\xrightarrow{\,({a_2,b_2})\,}\ 
                         \cdots\ 
                       \xrightarrow{\,({a_p,b_p})\,}v_p=\eg\,,
 \]
 satisfying 
 \begin{enumerate}
 \item[(1)] $b_1\ge b_2\ge\dotsb\ge b_p\,$, and
 \item[(2)] $a_i\ne a_j$ if $i\ne j\,$.
 \end{enumerate}
 This formula is multiplicity free as there is at most one such chain between any two
 permutations. 
\end{thm}

The Pieri formula was first stated by Lascoux and Sch\"utzenberger in~\cite{lasps}, where they suggested a proof by induction on both $p$ and $k$; this was later carried out in~\cite{manfsp} (but with a subtle error). The first proof that appeared was the geometric one in~\cite{sotprf}. A different formulation and approach can be found in~\cite{posqvp}. We present an inductive proof of the Pieri formula in $K$-theory, and thus correct the proof in~\cite{manfsp}.

\subsection{The Pieri-type formula}

\begin{dfn}\label{D:Pieri_Chain}
 {\rm 
 A {\em marked chain} in the $k$-Bruhat order is a saturated chain $\gamma$
  \begin{equation}\label{E:Marked_Chain}
     v=v_0\xrightarrow{\,({a_1,b_1})\,}
                   v_1\xrightarrow{\,({a_2,b_2})\,}\ 
               \cdots\ 
                      \xrightarrow{\,({a_q,b_q})\,}v_q=\eg\,,
         \;\;\;q=\ell(\gamma)\,,
  \end{equation}
 in the $k$-Bruhat order with some covers marked,
 which we often indicate by underlining their labels:
 $v_{i-1}\xrightarrow{\,\underline{(a_i,b_i)}\,}v_i$.

 A {\it Pieri chain} in the $k$-Bruhat order is a marked chain in the $k$-Bruhat
 order which satisfies the following four conditions.
 \begin{enumerate}
  \item[(P1)] $b_1\ge b_2\ge\dotsb\ge b_q\,$.\vspace{-2pt}
  \item[(P2)] If the $i$th cover
              $v_{i-1}\xrightarrow{\,\underline{(a_i,b_i)}\,}v_i$ 
              is marked, then $a_j\ne a_i$ for $j<i\,$.
  \item[(P3)] If the $i$th cover
              $v_{i-1}\xrightarrow{\,(a_i,b_i)\,}v_i$ 
              is not marked and $i+1\le q$, then $(a_i,b_i)\prec(a_{i+1},b_{i+1})$.
  \item[(P4)] If $b_1=\dotsb=b_r$ and $a_1>\dotsb>a_r$ for some 
              $r\ge 1$, then $(a_r,b_r)$ is marked.
 \end{enumerate}
}
\end{dfn}

For example, here is a Pieri chain in the $4$-Bruhat order on $S_7$
\[
 \Magenta{4}\Blue{261\,73}\Magenta{5}  
 \xrightarrow{\,\underline{(1,7)}\,}\ 
 \Blue{5}\Magenta{2}\Blue{61\,7}\Magenta{3}\Blue{4}
 \xrightarrow{\,(2,6)\,}\ 
 \Blue{536}\Magenta{1}\,\Blue{7}\Magenta{2}\Blue{4}
 \xrightarrow{\,(4,6)\,}\ 
 \Blue{53}\Magenta{6}\Blue{2}\,\Magenta{7}\Blue{14}
 \xrightarrow{\,\underline{(3,5)}\,}\ 
 \Blue{5}\Magenta{3}\Blue{72}\,\Magenta{6}\Blue{14}
 \xrightarrow{\,(2,5)\,}\ 
 \Blue{5672314}\,.
\]

We write $\fg_\lambda(x_1,\ldots,x_k)$ for the Grothendieck polynomial
$\fg_{v(\lambda,k)}(x)$, as it is symmetric in $x_1,\ldots,x_k$. 
Hence we have 
 \[
   \fg_{(1^p)}(x_1,\ldots,x_k)\ :=\ \fg_{c_{[k,p]}}(x)\qquad\mbox{and}\qquad
    \spa\spa\fg_{(p)}(x_1,\ldots,x_k)\ :=\ \fg_{r_{[k,p]}}(x)\,.
 \]
The $K$-theory Schubert classes represented by
$\fg_{v(\lambda,k)}(x)$ are pulled back from Grassmannian projections. 

\begin{thm}\label{kthm} 
  We have that
 \begin{equation}\label{kpieri}
   \fg_v(x)\,\fg_{(1^p)}(x_1,\ldots,x_k)=
   \sum(-1)^{\ell(\gamma)-p}\,\fg_{\eg}(x)\,,
 \end{equation}
 where the sum is over all Pieri chains $\gamma$ in the $k$-Bruhat order (on $S_\infty$) 
 that begin at $v$ and have $p$ marks.
 This formula has no cancellations.
\end{thm}

\begin{rem} Let us give a combinatorial argument for why the sum in (\ref{kpieri}) has a finite number of terms. Fix the permutation $v$ in $S_\infty$, and let $n$ be the smallest integer such that $n\ge k$ and $v(i)=i$ for $i>n$. The cover condition implies that the first cover $v=v_0\xrightarrow{\,({a_1,b_1})\,} v_1$ in a Pieri chain $\gamma$ (in the $k$-Bruhat order) is such that $b_1\le n+1$. Using the notation $(\ref{E:Marked_Chain})$ for $\gamma$, we deduce from condition (P1) that both $a_i$ and $b_i$, for $i=1,\ldots,q$, are bounded. On the other hand, we cannot have $(a_i,b_i)=(a_j,b_j)$ for $i\ne j$; indeed, using condition (P1) again, the above fact would imply that $b_l=b_i=b_j$ for all $l$ between $i$ and $j$, which would contradict the fact that $\gamma$ is increasing in Bruhat order. Hence $\gamma$ is given by a word with non-repeated letters $(a_i,b_i)$ (marked or not) in a finite alphabet. Clearly, there are a finite number of such words. 
\end{rem}

We interpret the conditions (P1)--(P4) for Pieri chains.
Condition (P1) is shared by both the Pieri-type formula for Schubert
polynomials (Condition (1) in Theorem~\ref{schub}) and the Monk formula for
Grothendieck polynomials (Condition~\eqref{E:inc_chain} of
Theorem~\ref{kmonk} implies that $b_1\ge b_2\ge \dotsb\ge b_p$).
The $p$ marked covers correspond to the $p$ covers in the Pieri-type formula for
Schubert polynomials, and condition (P2) is an analog to Condition (2) in
Theorem~\ref{schub}. 
For (P3), the unmarked covers behave like those in the Monk formula, analogous 
to Condition~\eqref{E:inc_chain} in Theorem~\ref{kmonk}.
For the last condition (P4), note that, by (P3), the first $r{-}1$ covers are forced
to be marked.
Thus (P4) states that the first cover that is not forced to be
marked by the previous conditions must be marked.

\begin{rem}\label{condp0}
 Consider a saturated chain $\gamma$ in the $k$-Bruhat
 order~\eqref{E:Marked_Chain} which admits a marking satisfying Conditions
 (P1)--(P4) for some number $p>0$ of marks.
 This can happen if and only if $\gamma$ satisfies Condition (P1) and the
 condition
\begin{enumerate}
\item[(P0)] For $i=2,\dotsc,\ell(\gamma)-1$, if $a_j=a_i$ for some $j<i$, then
       $(a_i,b_i)\prec(a_{i+1},b_{i+1})$.
\end{enumerate}
\end{rem}

In Section~\ref{S:uniqueness} we prove Theorem~\ref{T:unique}, which states that
there is at most one chain in the $k$-Bruhat order between two permutations
satisfying Conditions (P0) and (P1).
A consequence of this uniqueness is a different version of the Pieri-type
formula of Theorem~\ref{kthm}. Let us note that these results were not conjectured before; the only known special cases correspond to $p=1$ (the Monk formula \cite{lenktv} of Theorem \ref{kmonk}) and to $v$ being a Grassmannian permutation (cf. \cite{lencak} and Theorem \ref{pieri1} below). 

Write $v\xrightarrow{\,c(k)\,}w$ when there is a chain $\gamma$ from $v$ to $w$
in the $k$-Bruhat order satisfying conditions (P0) and (P1).
If $v\xrightarrow{\,c(k)\,}w$ with chain $\gamma$, then some covers in $\gamma$
are forced to be marked by conditions (P2)--(P4), while other covers are
prohibited from being marked.
Let $f(v,w)$ be the number of covers in $\gamma$ forced to be marked and
$p(v,w)$ the number of covers prohibited from being marked.
Set the binomial coefficient $\binom{n}{k}=0$ if $0\le k\le n$ does not hold.

\begin{cor}\label{compressedformula}  
  We have that
 \[
   \fg_v(x)\,\fg_{(1^p)}(x_1,\ldots,x_k)\ =\ 
   \sum(-1)^{\ell(w)-\ell(v)-p}\,
    \binom{\ell(w)-\ell(v)-f(v,w)-p(v,w)}{p-f(v,w)}\,\fg_w(x)\,,
 \]
  the sum over all permutations $w$ (in $S_\infty$) such that $v\xrightarrow{\,c(k)\,}w$.  
\end{cor}

\section{Some finer aspects of the Bruhat order}\label{S:Bruhat_Order}

We collect some results on the Bruhat order.
in Section~\ref{S:uniqueness}, we show that there is at most one chain in the
$k$-Bruhat order satisfying Conditions (P0) and (P1).
Our proof of the Pieri-type formula requires several rather technical results on
chains in the Bruhat order, which we give in Section~\ref{S:forbidden}.

We recall a characterization of the $k$-Bruhat order~\cite[Theorem A]{basspb}.
Given permutations $v,w\in S_n$, we have $v\leq_k w$ if and only if
 \begin{equation}\label{Characterization}
  \begin{minipage}{400pt}
   \begin{enumerate}
     \item $a\leq k<b$ implies that $v(a)\leq w(a)$ and $v(b)\geq w(b)$.
     \item $a<b$ with $v(a)<v(b)$ and $w(a)>w(b)$ implies that $a\leq k<b$. 
   \end{enumerate}
  \end{minipage}
 \end{equation}

\subsection{Uniqueness of Pieri chains}\label{S:uniqueness}

\begin{thm}\label{T:unique}
  Let $v$ and $w$ be permutations.
  There is at most one saturated chain in $k$-Bruhat order from $v$ to $w$
  satisfying Conditions {\rm (P0)} and {\rm (P1)}.  
\end{thm}

The weaker result that chains $\gamma$ satisfying  
$(a_i,b_i)\prec(a_{i+1},b_{i+1})$ for  $i=2,\dotsc,\ell(\gamma)-1$ are unique
(uniqueness of chains in the Monk formula in Theorem~\ref{kmonk}) was
proved in \cite{lenktv}.   
The special case when $\ell(\gamma)=p$ is part of Theorem \ref{schub},
for the Pieri formula in cohomology.  

\begin{proof}
 If an integer $i$ does not occur in an expression $x$, we {\it flatten} $x$
 by replacing each occurrence of $j$ with $j>i$ in $x$ by $j{-}1$.
 Let $\gamma$ be a chain in the $k$-Bruhat order from $v$ to $w$ in
 $S_n$ satisfying Conditions (P0) and (P1). 
 We assume without loss of generality that $v(i)\neq w(i)$ for all
 $i=1,\dotsc,n$.
 Indeed, if $v(i)=w(i)$, then $\gamma$ has no transposition involving position
 $i$. 
 If we restrict $v$ and $w$ to positions $i$ with $v(i)\neq w(i)$, flatten the
 results, and likewise flatten $\gamma$, then we obtain a chain $\gamma'$ from
 $v'$ to $w'$ in $S_{n'}$ which satisfies Conditions (P0) and (P1),
 and where $v'(i)\neq w'(i)$ for all $i=1,\dotsc,n'$.
 We recover $v,w$, and $\gamma$ from $v',w'$, $\gamma'$, and the set 
 $\{i\mid v(i)=w(i)\}$.
 
 If $v(i)\neq w(i)$ for all $i=1,\dotsc,n$, then Condition (P1) implies that
 $\gamma$ has the form 
\[
  (a_1,n),\dotsc,(a_r,n),\, (b_1,n{-}1),\dotsc,(b_s,n{-}1),\,\dotsc,
  (c_1,k{+}1),\dotsc,(c_t,k{+}1)\,,
\]
 where $r,s,t>0$. Observe that $a_r=v^{-1}(w(n))$.
 We will show that $a_1,\dotsc,a_{r-1}$ are also determined by $v$ and $w$, which
 will prove the theorem by induction on $n$, as any final segment of $\gamma$ is
 a chain satisfying Conditions (P0) and (P1).
 More precisely, set $\alpha:=v(a_r)=w(n)$, $\beta:=v(n)$, and $a_0:=n$. 
 Then we will show that whenever $i\not\in\{a_0,\dotsc,a_r\}$, then
 either $v(i)<\alpha$ or $v(i)>\beta$.
 This will imply that $\alpha=\beta-r$ and thus $v(a_i)=\beta-i$, for
 $i=0,\dotsc,r$.

 Since $\gamma$ is a chain in the Bruhat order, 
 $\alpha=v(a_r)<\dotsb<v(a_1)<v(a_0)=\beta$ and the transposition $(a_j,n)$ in
 $\gamma$ interchanges $v(a_j)$ with $v(a_{j-1})$, which is in position $n$. 
 Suppose by way of contradiction that there is some $i\not\in\{a_0,\dotsc,a_r\}$
 with $\alpha<v(i)<\beta$.
 Let $j$ be the index with $v(a_j)<v(i)<v(a_{j-1})$.
 By the cover condition applied to the cover $(a_j,n)$, we have $i<a_j$.
 At this point in the chain $\gamma$, the value in position $i$ is less than
 the value in position $a_j$, and so by the characterization of the $k$-Bruhat
 order~\eqref{Characterization}, this remains true for subsequent permutations in the chain.

 Let $(i,m)$ be the first transposition in $\gamma$ which involves position $i$
 and $u$ the permutation to which this transposition is applied. 
 Then $u(i)<u(a_j)$, as we observed.
 Since $i<a_j<m$, the cover condition implies that $u(m)<u(a_j)$.

 As $v\leq_k u$ and $k<m$, we have
\[
   v(m)\ \geq\ u(m)\ >\ u(i)\ =\ v(i)\ >\ v(a_j)\,.
\]
 Thus somewhere in the chain $\gamma$ after $(a_j,n)$ but before $(i,m)$, the
 relative order of the values in positions $a_j$ and $m$ is reversed.
 By Condition (P1) and the cover condition, this happens by applying a
 transposition $(l,m)$ with $a_j\leq l$.

 Since $i<l$, there are two transpositions $(l',m),(l'',m)$ which are adjacent
 in the chain $\gamma$, occur weakly between $(l,m)$ and $(i,m)$, and 
 satisfy $l''<l\leq l'$.
 Assume that this is the first such pair and let $x$ be the permutation to which
 $(l',m)$ is applied. 
 By Condition (P0), the transposition $(l',m)$ is the first occurrence
 of $l'$ in the chain $\gamma$,
 so $a_j\ne l'$ (which means that $a_j<l'$), as $a_j$ occured 
 before in the transposition $(a_j,n)$.  
 We claim that $x(l')<x(a_j)$.
 Indeed, if $l=l'$, then $x(l')<x(a_j)<x(m)$, by the definition of $l$,
 and if $l'\neq l$, then $x(l')<x(m)<x(a_j)$.

 Since the transposition $(l',m)$ is the first occurrence
 of $l'$ in the chain $\gamma$, we have $x(l')=v(l')$.
 Thus, $v(l')$ becomes the value in position $m$ after applying $(l',m)$, so we have 
\[
     v(l')\ \geq\ u(m)\ >\ u(i)\ =\ v(i)\ >\ v(a_j)\,.
   \]
 Since $a_j<l'$ and $v\leq_k x$,~\eqref{Characterization} implies that
    $x(a_j)<x(l')$, a contradiction.
 \end{proof}

\begin{rem}
 If the permutation $v$ has no ascents in positions $1,\ldots,k-1$, 
 then if $v\xrightarrow{\,c(k)\,}w$, the chain
 $\gamma$ from $v$ to $w$ is a Monk chain.
 Indeed, we cannot have $(a_i,b_i)\not\prec(a_{i+1},b_{i+1})$ in $\gamma$ for
 then $b_{i+1}=b_i$ and $a_{i+1}<a_i$.
 But the value in position $a_i$ is $v(a_i)$, by Condition (P0), and this 
 is less than the value in position $a_{i+1}$, by the special form of $v$ and
 the characterization of the $k$-Bruhat order~\eqref{Characterization}.
 Hence, the permutations indexing the Grothendieck polynomials with nonzero
 coefficients in the expansion of $\fg_v(x)\,\fg_{(1^p)}(x_1,\ldots,x_k)$ are
 precisely the permutations relevant to the expansion of
 $\fg_v(x)\,\fg_{s_k}(x)$ (the Monk-type formula in Theorem \ref{kmonk}) which
 differ from $v$ in at least $p$ positions from $1$ to $k$.  

 If the permutation $v$ has no descents in positions
 $1,\ldots,k-1$, then 
 a stronger version of  Theorem \ref{T:unique} holds, in which
 Condition (P0) is dropped from the hypothesis. 
 Indeed, there is always at most one position from 1 to $k$ with which we can
 transpose a given position greater than $k$ such that the cover condition
 holds. 
 In general, not all chains obtained in this case satisfy Condition (P0). 
 We also note that the statement in Theorem \ref{T:unique} is false
 if Condition (P0) is dropped from the hypothesis. 
\end{rem}

The proof of Theorem \ref{T:unique} provides the following algorithm to decide
if $v\xrightarrow{\,c(k)\,}w$ and if so, produces the chain $\gamma$.

\begin{alg}\hfill
 \begin{itemize}
  \item[{\rm Step 1.}] Let $m:=n$ and $\gamma:=\emptyset$.

  \item[{\rm Step 2.}] If $v(m)=w(m)$ then $A_m:=\emptyset$; go to {\rm Step 6}.

  \item[{\rm Step 3.}] Let 
\[
  A_m\ :=\ 
   \{i \mid i\le k\,,\spa v(i)\ne w(i)\,,\spandsp w(m)\le v(i)< v(m)\}\,.
\]

  \item[{\rm Step 4.}] Write the elements in $A_m$ as $\{a_1,\dotsc,a_r\}$, where
 $v(m)>v(a_1)>\dotsb>v(a_r)$, and then set
\[
     v\ :=\ v (a_1,m)(a_2,m)\dotsb(a_r,m)\qquad\mbox{\rm and}\qquad
     \gamma\ :=\ \gamma | (a_1,m),(a_2,m),\dotsc,(a_r,m)\,.
\]
  (Here, $|$ means concatenation.)

  \item[{\rm Step 5.}]  If $v(m)\ne w(m)$, or any multiplication by a transposition
 $(a_i,m)$ in Step {\rm 4} violates Condition {\rm (P0)}, then output
 ``no such chain''. STOP. 

  \item[{\rm Step 6.}] Let $m:=m-1$; if $m>k$ then go to Step {\rm 2}.

  \item[{\rm Step 7.}] If $v=w$ then output the chain $\gamma$
    else output ``no such chain''. STOP.  
 \end{itemize}
\end{alg} 

\begin{rem}
    The branching rule in the tree of saturated chains (in $k$-Bruhat
    order) satisfying Conditions (P0) and (P1) which start at a given
    permutation is simple. Indeed, if we are at the beginning of the
    chain, any transposition $(a,b)$ with $a\le k<b$ that satisfies the cover
    condition can be applied. Otherwise, assuming $(c,d)$ is the previous
    transposition, the current transposition $(a,b)$ has to
    satisfy the following extra conditions: (1) $b\le d$; (2) if $b=d$ and there
    is a transposition $(c,\cdot)$ before $(c,d)$, then $a>c$.  
\end{rem}

\begin{ex}
 Consider multiplying $\fg_{21543}(x)$ by $\fg_{(1^2)}(x_1,x_2,x_3)$. 
 There is a unique chain in the 3-Bruhat order from $215436$ to $426315$ which
 satisfies Condition (P1).
\[
  \Blue{21}\Magenta{5}\,\Blue{43}\Magenta{6}\ 
   \xrightarrow{\,(3,6)\,}\ 
  \Magenta{2}\Blue{16\,4}\Magenta{3}\Blue{5}\ 
   \xrightarrow{\,(1,5)\,}\ 
  \Blue{3}\Magenta{1}\Blue{6\,4}\Magenta{2}\Blue{5}\ 
   \xrightarrow{\,(2,5)\,}\ 
  \Magenta{3}\Blue{26}\,\Magenta{4}\Blue{15}\ 
   \xrightarrow{\,(1,4)\,}\ 
  \Blue{426315}\ .
\]
 This chain has two markings that satisfy conditions (P2)--(P4) for $p=2$.
\[
  \Bigl(\underline{(3,6)},\underline{(1,5)},(2,5),(1,4)\Bigr)
     \qquad\mbox{and}\qquad
  \Bigl(\underline{(3,6)},(1,5),\underline{(2,5)},(1,4)\Bigr)
\]
 Hence, the coefficient of $\fg_{426315}(x)$ in the product 
 $\fg_{21543}(x)\cdot \fg_{(1^2)}(x_1,x_2,x_3)$ is 2.
\end{ex}

\subsection{Forbidden subsequences in chains in the $k$-Bruhat
             order}\label{S:forbidden} 
The Bruhat order is Eulerian~\cite{Ve71}, so 
every interval of length two has 2 maximal chains.
This defines a pairing on chains of length two.
For convenience, we may represent a chain in the Bruhat order by the sequence of
its covering transpositions; thus the chain 
\[
   v\ \xrightarrow{\,(a,b)\,}\ 
   v'\ \xrightarrow{\,(c,d)\,}\ v''\quad\mbox{may be written}\quad
   \bigl( (a,b),(c,d) \bigr)\,.
\]

\begin{lem}\label{L:length2}
 Let $\gamma$ be a chain of length two in the Bruhat order from $v$ to $w$.
 \begin{enumerate}
  \item If $\gamma=\big((a,b),(c,d)\big)$ with
    $a,b,c,d$ distinct numbers, then $\big((c,d),(a,b)\big)$ is
    another chain from $v$ to $w$. 
  \item Suppose that $j<k<l$.
   \begin{enumerate}
    \item If $\gamma=\bigl( (k,l),(j,l) \bigr)$, then the other chain is
                    $\bigl( (j,k),(k,l) \bigr)$.
    \item If $\gamma=\bigl( (j,l),(j,k) \bigr)$, then the other chain is
                    $\bigl( (j,k),(k,l) \bigr)$.
    \item If $\gamma=\bigl( (j,l),(k,l) \bigr)$, then the other chain is
                    $\bigl( (k,l),(j,k) \bigr)$.
    \item If $\gamma=\bigl( (j,k),(j,l) \bigr)$, then the other chain is
                    $\bigl( (k,l),(j,k) \bigr)$.
   \end{enumerate}
 \end{enumerate}
\end{lem}

Transformations involving (1) will be called {\em commutations} and
those involving (2) {\em intertwining relations}.
Note that the second chains in (2)(a) and (2)(b) are the same, as are the second
chains  (2)(c) and (2)(d).
Thus chains of either form $\bigl( (j,k),(k,l)\bigr)$ or
$\bigl( (k,l),(j,k)\bigr)$ intertwine in one of two distinct ways. 

\begin{rem}\label{R:M1}
 The omission in the proof of the cohomology Pieri formula
 in~\cite{manfsp} arises from intertwining transpositions and it occurs in the
 displayed formula at the bottom of page~94 (which also contains a small
 typographic omission---$\zeta\zeta'$ should be replaced by
 $\zeta\zeta't_{mq}$). 
 There, a transposition $(a,b)$ is written $t_{ab}$.
 In that formula, transpositions $t_{i_rm}t_{mq}$ with $i_r<m<q$ are intertwined
 as in Lemma~\ref{L:length2}2(b), and the case when they intertwine as in 2(a)
 is neglected.
 This neglected case does occur when (in the notation of~\cite{manfsp})
 $w=32154$, $m=p=4$, $q=5$, and $u\in S_{3,3}(w)$ is 
 $45213 = wt_{24} t_{34} t_{15}$.
 We discuss this further in Remarks~\ref{R:M2} and~\ref{R:M3}.
\end{rem}

We give four lemmas prohibiting certain chains in the $k$-Bruhat order.

\begin{lem}\label{L:Prohibition_1}
  A saturated chain in the $k$-Bruhat order satisfying {\rm (P0)} and 
  {\rm (P1)} cannot contain a subsequence of the form
  $(j,m),\dotsc,(i,m),\dotsc,(i,l)$ with $i<j\leq k<l<m$.
\end{lem}

\begin{proof}
 We show that if a saturated chain $\gamma$ in the $k$-Bruhat order satisfying
 {\rm (P0)} and {\rm (P1)} contains a subsequence   
 $(j,m),\!$ $\dotsc,\!$ $(i,m),\!$ $\dotsc,\!$ $(i,l)$ with 
 $i<j\leq k<l<m$, then there is a saturated chain in the $k$-Bruhat order
 $(j',m')(i',m')(i',l')$ with $i'<j'\leq k<l'<m'$.
 If this chain begins at a permutation $v$, then the cover condition
 implies the contradictory inequalities 
 $v(l')<v(m')$ and $v(m')<v(l')$.

 Removing an initial segment from $\gamma$, we may assume that it begins with
 $(j,m)$ and no transposition $(h,m)$ between $(j,m)$ and $(i,m)$  
 satisfies $i<h$.
 Removing a final segment, we may assume that it ends with $(i,l)$
 and its only transpositions involving $i$ are $(i,m)$ and $(i,l)$.
 This gives a chain which we still call $\gamma$ that satisfies Conditions
 {\rm (P0)} and {\rm (P1)}.

 If $(j,m)$ is not next to $(i,m)$, then it is next to a
 transposition $(h,m)$ with $h<i<j$.
 By Lemma~\ref{L:length2}(2)(a), we may replace
 $(j,m)(h,m)$ by $(h,j)(j,m)$ and then remove $(h,j)$, obtaining another 
 chain satisfying {\rm (P0)} and {\rm (P1)}. 
 Continuing in this fashion gives a chain satisfying {\rm (P0)} and {\rm (P1)}
 in which $(j,m)$ is adjacent to $(i,m)$.

 We reduce to the case $h_1>h_2>\dotsb>h_r>i$, where 
 $(h_1,l),(h_2,l),\dotsc,(h_r,l),(i,l)$ is the subchain of transpositions
 involving $l$.
 While there is a transposition $(a,l)$ followed by $(b,l)$ with
 $a<b$, pick the rightmost such and use the intertwining relation of
 Lemma~\ref{L:length2}(2)(c) to replace $(a,l)(b,l)$ by $(b,l)(a,b)$.
 Next, commute $(a,b)$ to the end of the chain and remove it.
 This is possible because the positions $h_1,\ldots,h_r,i$ are all distinct.

 Since $h_1>h_2>\dotsb>h_r>i$ and the chain ends with these
 transpositions involving $l$, Condition (P0) implies that no
 $h_i$ occurs elsewhere in the chain.
 Thus all transpositions between $(i,m)$ and $(h_1,l)$ can be
 commuted to the end of the chain and removed, while the transpositions 
 $(h_1,l),(h_2,l),\dotsc,(h_r,l)$ can be commuted to the front of
 the chain and removed.
 This gives a chain in the Bruhat order of the form
 $(j',m'),(i',m'),(i',l')$ with $i'<j'<l'<m'$, which, as we observed, is
 forbidden.
\end{proof}

\begin{lem}\label{L:Prohibition_2}
  Suppose that $\gamma$ is a saturated chain in the $k$-Bruhat order from $v$ to
  $u$ satisfying Condition {\rm (P1)}, 
  and that there is a saturated chain above $u$ in the $k$-Bruhat order of the form
 \begin{equation}\label{E:l-segment}
    u\xrightarrow{\,({i_1,l})\,}
                   u_1\xrightarrow{\,({i_2,l})\,}\ 
               \cdots\ 
                      \xrightarrow{\,({i_r,l})\,}u_r\,.
 \end{equation}
 Then the concatenation of $\gamma$ with this chain cannot contain a segment
 of the form
 \begin{equation}\label{E:Forbidden_Segment}
     (i,m),\,\dotsc,\,(i,l),(j,l)\,,
 \end{equation}
 where $i<j\leq k<m\leq l$, and no other transposition
 in~\eqref{E:Forbidden_Segment} involves  $i$.  
\end{lem}

\begin{proof}
 We may assume that the concatenation begins with $(i,m)$,
 ends with $(j,l)$, and that no transposition between $(i,m)$ and $(i,l)$
 involves $i$. 
 Let $w$ be the permutation to which $(i,l)$ is applied.
 Then the covers $v\lessdot v(i,m)$ and $w\lessdot w(i,l)$, together with the
 hypothesis on $i$ imply that 
\[
   v(i)\ <\ v(m)\ =\ w(i)\ <\ w(l)\ \leq\ v(l)\,.
\]
 In particular, $l\neq m$ and so $i<j<m<l$.
 This implies that $(i,m)$ belongs to the chain $\gamma$ and $(i,l),(j,l)$
 to the chain~\eqref{E:l-segment}, by Condition (P1).

 We claim that there is a number $y$ in one of the positions $j{+}1,\dotsc,m$ in
 $w$ that satisfies $v(j)<y<v(m)$.
 This will lead to a contradiction to the cover condition for the cover $(j,l)$.
 Indeed, assumptions on $i$ and $j$ imply that the numbers exchanged at the
 cover $(j,l)$ are $v(j)<v(m)$ with $v(j)$ in position $j$ and
 $v(m)$ in position $l$.

 We prove the claim by induction on the length of the chain.
 After the first cover in the chain, $v(i)$ occupies position $m$.
 Since we have $v(j),v(i)<v(m)$ and $i<j<m$, the cover condition implies that
 $v(j)<v(i)$.

 Suppose that a subsequent permutation in the chain has a value $y$ in a
 position $h$ with $j\leq h\leq m$ and $v(j)<y<v(m)$.
 If $k<h$, then any transposition $(n,h)$ that interchanges $y$ with a number
 smaller than $v(j)$ must have $j<n$ by the cover condition, and so $y$ remains
 at a position $n$ with $j<n\leq m$.
 Let us now turn to the case $h\leq k$. Assume that the next cover in the
 chain is $(h,n)$.
 If this cover comes from the chain $\gamma$, then $k<n\leq m$, as $\gamma$
 satisfies (P1), and so $y$ remains at a position $n$ with $j<n\leq m$.
 If the cover comes from the segment~\eqref{E:l-segment}, then it places $y$ in
 position $l$.
 This implies that subsequent permutations in the chain have values at most $y$
 in position $l$, which contradicts $v(m)$ occupying that position in the
 penultimate step in the chain.
 This proves the claim and completes the proof of the lemma.
\end{proof}

\begin{lem}\label{L:technical}
 Let $\gamma$ be a saturated chain in the $(k{-}1)$-Bruhat order from $v$ to
 $u$ satisfying Conditions {\rm (P0)} and {\rm (P1)}.
 Assume that $\gamma$ has the form
 \begin{equation}\label{E:string}
   \gamma=\gamma'\,|\,(i_0,k),(i_1,k),\dotsc,(i_r,k) \,,
 \end{equation}
 and, for $t=1,\dotsc,r$, the transpositions $(i_t,k),(k,l)$ intertwine as
 in Lemma~${\rm \ref{L:length2}}(2)({\rm a})$. 
 Also assume that any transposition $(i_r,m)$ in $\gamma'$ has $m>l$. 
 Then any transposition $(i_t,m)$ in
 $\gamma'$ has $m>l$ for $t=0,\dotsc,r$.  
\end{lem}

\begin{proof}
 We use induction on $r$ to show that any transposition $(i_0,m)$ in $\gamma'$
 has $m>l$. 
 The assumption on $\gamma'$ is the case $r=0$.
 Assume $r\ge 1$, and consider the concatenation $\gamma|(k,l)$ of $\gamma$
 with $(k,l)$. 
 Let us intertwine $(i_t,k)$ with $(k,l)$, for $t=r,r-1,\ldots,0$. 
 If $(i_0,k)$ intertwines with $(k,l)$ as in
 Lemma~\ref{L:length2}(2)(a), then we obtain the chain
 \begin{equation}\label{form1}
   \gamma'\,|\, (k,l),(i_0,l),\,\dotsc,\,(i_r,l)\,.
 \end{equation}
 Otherwise, 
 $(i_0,k)(k,l)$ becomes $(i_0,l)(i_0,k)$. 
 We then commute $(i_0,k)$ to the right and remove it to obtain the chain
 \begin{equation}\label{form2}
   \gamma'\,|\, (i_0,l),\,\dotsc,\,(i_r,l)\,.
 \end{equation}
If $\gamma'$ contains a transposition $(i_0,m)$, then we must have $i_0<i_1$, by
Condition (P0). 
Given that $\gamma'$ satisfies Conditions (P0) and (P1), we conclude that
$m>l$. 
Indeed, otherwise we have $i_0<i_1\le k<m\le l$, so the segment
\[ 
   (i_0,m),\,\dotsc,\,(i_0,l),(i_1,l)\,,
\]
in~\eqref{form1} or~\eqref{form2} has the form forbidden by 
Lemma~\ref{L:Prohibition_2}, unless there is a transposition $(i_1,n)$ between
$(i_0,m)$ and $(i_0,l)$. 
But then the induction hypothesis
applies, and we have $m\ge n>l$, by Condition (P1).  
\end{proof}

\begin{lem}\label{L:quad}
  A chain in the $(k{-}1)$-Bruhat order satisfying {\rm (P0)} and {\rm (P1)} 
  cannot contain a segment of the form
 \begin{equation}\label{E:quad}
     (i,l),\,\dotsc,\,(h,l),\,\dotsc,\,(h,m),\,\dotsc,\,(i,k)\,,
 \end{equation}
 where there are no other transpositions between $(i,l)$ and $(i,k)$ which
 involve $i$.
\end{lem}

\begin{proof}
 Suppose that $\gamma$ is such a chain.
 Truncate $\gamma$ to begin at a permutation $v$ with the cover 
 $v\lessdot v(i,l)$ and to end with the step $z\lessdot z(i,k)$.
 Furthermore, the assumptions on $i$ imply that $z(i)=v(l)$.
 The covers  $(i,l)$ and $(h,l)$ imply that
\[
   v(h)\ <\ v(i)\ <\ v(l)\,.
\]

 Since $\gamma$ satisfies Condition (P1), we have $i\neq h$, $m\neq l$, and 
 $k\leq m< l$.
 We must have $i<h$ for if $i>h$, then 
 $(i,l),\dotsc,(h,l),\dotsc,(h,m)$ is forbidden by Lemma~\ref{L:Prohibition_1}. 

 We next argue that $m\neq k$.
 Suppose otherwise that $m=k$, so $(h,m)$ in~\eqref{E:quad} becomes
 $(h,k)$ and let $y\lessdot y(h,k)=:y'$ be the corresponding cover.
 By Condition (P0), there is a permutation $x$ after $y$ and two covers
\[
   x\ \xrightarrow{\,(h',k)\,}\ x'\ 
      \xrightarrow{\,(h'',k)\,}\ x''\,,
\]
 where $h''<h<h'$.
 By Condition (P0) and the cover condition, $v(h')=x(h')<x(k)$.
 Since $y'= y(h,k)$, we have $y'(k)<y'(h)$.
 Since $y'\leq_{k-1}x$, we have $x(k)\leq y'(k)$ and 
 $y'(h)\leq x(h)$.
 Putting these together gives
 \begin{equation}\label{Eq:String}
    x(h')\ <\ x(h)\,.
 \end{equation}
 Since $x'\leq_{k-1}z$, we have $x'(k)\geq z(k)$, and also 
 $v(h')=x(h')=x'(k)>v(l)$.
 Since $v(l)>v(h)$, we have $v(h)<v(h')$.
 However, $v\leq_{k-1}x$, and so the
 characterization~\eqref{Characterization} of $\leq_{k-1}$ implies
 that $x(h)<x(h')$, contradicting~\eqref{Eq:String}.
 
 Thus $i<h<k<m<l$.
 We determine the values of $v$ at these positions.
 Let $(h,m)$ be the last transposition in the chain $\gamma$
 involving $h$, with cover $y\lessdot y(h,m)=:y'$.
 By the cover condition for the covers $(i,l)$ and $(h,l)$, neither
 $v(k)$ nor $v(m)$ lies between $v(h)$ and $v(l)$.
 
 Recall that the last cover in the chain is $z\lessdot z(i,k)$ and 
 $v(l)=z(i)<z(k)$.
 Since $v\leq_{k-1}z$, we have $z(k)\leq v(k)$ and so $v(l)<v(k)$.

  As $v(h)<v(k)$ and $v\leq_k y'$, the
 characterization~\eqref{Characterization} of $\leq_k$ implies that
 $y(m)=y'(h)<y'(k)=y(k)$.   
 But $v\leq_{k-1}y$, so the
 characterization~\eqref{Characterization} of 
 $\leq_{k-1}$ implies then $v(m)<v(k)$.
 We conclude that 
\[
    v(h)\ <\ v(i)\ <\ v(l)\ <\ v(m)\ <\ v(k)\,.
\]

 As $v\leq_{k-1}y$, we have $y(m)\leq v(m)$.
 We argue that it is impossible to have either 
 $v(l)<y(m)$ or $v(l)>y(m)$.
 (Since $i\neq m$ and $v(l)=y(i)$, we cannot have $v(l)=y(m)$.)
 
 Suppose that $v(l)<y(m)$.
 Since $(h,m)$ is the last transposition involving $h$, we have $z(h)=y(m)$.
 Since $z(i)=v(l)<z(h)$ and $z\lessdot z(i,k)$ is the last cover in the chain,
 the cover condition implies that $z(k)<z(h)$. 
 Since $y'(h)<y'(k)=v(k)$, there is some step
 $u\lessdot u(j,k)=:u'$ in the chain such that 
 $u(j)<u(h)=y(m)<u(k)$.

 The cover condition implies that $h<j$.
 As $u(h)>u(j)$ and $v\leq_{k-1}u$,~\eqref{Characterization} implies
 that $v(h)>v(j)$.
 Since $u(j)=u'(k)$ and $u'\leq_{k-1}z$, we have $u(j)\geq z(k)$.
 But then $u(j)>v(l)>v(h)>v(j)$, and so there had to be
 a previous transposition $(j,\cdot)$ in the chain.

 By Condition (P0), there is a permutation $x$ in the chain after $u$ and two
 covers
\[
   x\ \xrightarrow{\,(\jmath',k)\,}\ x'\ 
      \xrightarrow{\,(\jmath'',k)\,}\ x''\,,
\]
 with $\jmath''<j<\jmath'$.
 By Condition (P0) again, $v(\jmath')=x(\jmath')(=x'(k))$.
 Since $x'\leq_{k-1}z$, $z(k)>v(l)$, and $v(l)>v(h)$, we must have
 $v(\jmath')>v(h)$. 
 Then~\eqref{Characterization} implies that 
 $x(h)<x(\jmath')$.
 On the other hand, $u'(k)\geq x(k)$ and so $u(j)> x(\jmath')$
 as $u(j)=u'(k)$ and $x(\jmath')=x'(k)$.
 Since $x(h)=y(m)>u(j)$, we have found a contradiction to the assumption that
 $x(h)<x(\jmath')$. 

 Suppose now that $v(l)>y(m)$.
 Since $v(m)>v(l)$, there is some cover $t\lessdot t(a,m)=:t'$ in the chain
 before $y$ where
\[
  t'(m)=t(a)\ <\ v(l)\ <\ t(m)\,.
\]
 Since $v(l)=t(i)$, the cover condition implies that $i<a$.
 Since $t'(k)=v(k)$, we have
\[
   v(l)=t'(i)\ <\ t'(a)\ <\ t'(k)=v(k)\,.
\]

 Thus, if $w\lessdot w(b,k)$ is the first cover in the chain involving $k$, then
 $t'\leq_k w$ and~\eqref{Characterization} implies that
\[
   v(l)=w(i)\ <\ w(a)\ <\ w(k)=v(k)\,.
\]
 On the other hand, the cover condition on $z\lessdot z(i,k)$ implies that
\[
   v(l)=z(i)\ <\ z(k)\ <\ z(a)\,.
\]
 It follows that there is some cover $u\lessdot u(c,k)=:u'$ in the chain where
 $u(a)<u(k)$ but $u'(a)>u'(k)$.
 We will rule out both cases of $c\neq a$ and $c=a$.

 If $c=a$, then by Condition (P0) there is a permutation $x$ in the chain after
 $u$ and two covers 
\[
   x\ \xrightarrow{\,(a',k)\,}\ x'\ 
      \xrightarrow{\,(a'',k)\,}\ x''\,,
\]
 with $a''<a<a'$.
 By Condition (P0), $x(a')=u(a')=v(a')$.
 Note that  as $u'\leq_{k-1}x$, 
\[
   u(a')\ =\ x(a')\ <\ x(k)\ \leq\ u'(k)\ =\ u(a)\,.
\]
 Since $v\leq_{k-1}u$ and $a<a'$,~\eqref{Characterization} implies
 that $v(a')<v(a)$. 
 But this is a contradiction, as we had $v(a)\leq t(a)<v(l)$ and 
 $v(l)<z(k)\leq x'(k)=x(a')=v(a')$.

 We are left with the case $c\neq a$.
 Then there is no cover $(a,k)$ in the chain and so
 $w(a)=z(a)$.
 By our definition of $c$, we have
\[
   u(c)\ <\ u(a)\ <\ u(k)\,,
\]
 and the cover condition forces $a<c$.
 Since $t\leq_{k-1}u$,~\eqref{Characterization} implies that $t(c)<t(a)$.
 However, $t(a)<v(l)$ and also $v(l)<z(k)\leq u'(k)=u(c)$ and so
 we conclude that $t(c)\neq u(c)$ and so there had to be a previous
 transposition $(c,\cdot)$ in the chain.

 By Condition (P0), there is a permutation $x$ in the chain after
 $u$ and two covers 
\[
   x\ \xrightarrow{\,(c',k)\,}\ x'\ 
      \xrightarrow{\,(c'',k)\,}\ x''\,,
\]
 with $c''<c<c'$.
 By Condition (P0), this is the first transposition involving $c'$,
 and so
\[
  v(c')\ =\ u(c')\ =\ x(c')\ <\ x(k)\ \leq\  u'(k)\ =\ u(c)\,;
\]
 the inequalities are from the cover $x\lessdot x(c',k)$ and $u'\leq_{k-1} x$.
 Thus $a<c<c'$ and $u(a)>u(c)>u(c')$.
 Since $v\leq_{k-1} u$,~\eqref{Characterization} implies that
\[
  v(a)\ >\ v(c)\ >\ v(c')\,.
\]
 On the other hand, we have
\[
  v(c')\ =\ x(c')\ =\ x'(k)\geq z(k)>v(l)\,,
\]
 and $v(l)>t(a)\geq v(a)$, a contradiction.
\end{proof}

\section{Proof of the Pieri-type formula}\label{S:proof}

Fix a permutation $v$ throughout.
For $1\leq p\leq k$, let $\Gamma_{k,p}$ be the set of Pieri chains in the
$k$-Bruhat order which begin at $v$ and have $p$ marks.
For $\gamma\in\Gamma_{k,p}$, set $\sgn(\gamma)=(-1)^{\ell(\gamma)-p}$.
Let $\Gamma'_{k,p}$ be the set of chains which are the concatenation of a Pieri
chain $\pi$ in $\Gamma_{k-1,p}$ with a Monk chain $\mu$ satisfying the
conditions of Theorem~\ref{xkmonk} for the multiplication of 
$\fg_{{\rm end}(\pi)}(x)$ by $x_k$.
For $\gamma=\pi|\mu\in\Gamma'_{k,p}$, let
$\sgn(\gamma)=\sgn(\pi)\sigma(\mu)$, where $\sigma(\mu)$ is the sign of $\mu$ in
Theorem~\ref{xkmonk}. 
These parts $\pi$ and $\mu$ are called the {\it Pieri-} and 
{\it Monk- chains} of $\gamma$. 

We prove Theorem~\ref{kthm} by double induction on $0\leq p\leq k$, as explained below. (If $p>k$, then we have $\fg_{(1^p)}(x_1,\ldots,x_k)=0$.) The case $p=1$ is the Monk
formula~\cite{lenktv} of Theorem~\ref{kmonk}.  Assume that the formula of Theorem~\ref{kthm} holds for $1,\ldots,k-1$ (as values of $k$) and all $p$; the base case $k=1$ is given by the Monk formula, as discussed above, or simply by the formula of Theorem \ref{xkmonk}. The transition formula~\cite{lastgp,lenktv} for $\fg_{(1^p)}(x_1,\ldots,x_k)$ is 
 \begin{equation} \label{trans}
  \begin{array}{rcl}
     0&=&\fg_{(1^p)}(x_1,\dotsc,x_k)\ -\ \fg_{(1^p)}(x_1,\dotsc,x_{k-1})\\
       &&-x_k\fg_{(1^{p-1})}(x_1,\ldots,x_{k-1})
         +x_k\fg_{(1^p)}(x_1,\ldots,x_{k-1})\,. \rule{0pt}{14pt}
  \end{array}
 \end{equation}
Let us multiply this by $\fg_v(x)$. Based on our induction hypotheses and the formula of Theorem \ref{xkmonk}, we have 
 \begin{equation}\label{E:induction}
  \begin{array}{rcll}
   \fg_v(x)\cdot \fg_{(1^p)}(x_1,\ldots,x_{k-1})&=&
     {\displaystyle 
        \sum_{\gamma\in\Gamma_{k-1,p}} \sgn(\gamma)\fg_{\eg}(x)\,,}\\
   \fg_v(x)\cdot \fg_{(1^p)}(x_1,\ldots,x_{k-1})\cdot x_k&=&
    {\displaystyle 
         \sum_{\gamma\in\Gamma'_{k,p}} \sgn(\gamma)\fg_{\eg}(x)}\,,
       &\mbox{and}
          \rule{0pt}{15pt}
  \end{array}   
 \end{equation}
the analog of the last formula for $p{-}1$ instead of $p$. Hence, it suffices to show that 
 \begin{equation}\label{E:esum}
   \begin{array}{rcl}
   0&=&{\displaystyle 
        \sum_{\gamma\in\Gamma_{k,p}} \sgn(\gamma)\fg_{\eg}(x)\ -\ 
       \sum_{\gamma\in\Gamma_{k-1,p}} \sgn(\gamma)\fg_{\eg}}(x)\\
     && {\displaystyle 
          -\  \sum_{\gamma\in\Gamma'_{k,p-1}} \sgn(\gamma)\fg_{\eg}(x)\ +\ 
        \sum_{\gamma\in\Gamma'_{k,p}} \sgn(\gamma)\fg_{\eg}(x)\,.}
          \rule{0pt}{19pt}
   \end{array}
 \end{equation}
Indeed, based on the transition formula above, this would imply that 
\[
  \fg_v(x)\cdot\fg_{(1^p)}(x_1,\dotsc,x_k)\ =\ 
   \sum_{\gamma\in\Gamma_{k,p}}\sgn(\gamma)\fg_{\eg}(x)\,,
\]
which is the formula of Theorem~\ref{kthm}.
The claim that the sum on the right is without cancellations is immediate:
the sign of a term in that sum depends only upon the length of the permutation
with which it ends.
\medskip

We prove the formula~\eqref{E:esum} by giving a matching on chains in the
multiset 
 \begin{equation}\label{E:chain_indices}
  \Gamma\ :=\  \Gamma_{k,p}\ \cup\ \Gamma_{k-1,p}\ \cup\ 
   \Gamma'_{k,p}\ \cup\ \Gamma'_{k,p-1}\,,
 \end{equation}
which matches chains having the same endpoint but different signs.
This is a multiset, as $\Gamma_{k,p}$ and $\Gamma_{k-1,p}$ may have chains in
common. 
We actually define two different matchings on {\em subsets} $A$ and $B$ (as opposed to submultisets) of~\eqref{E:chain_indices} whose union is $\Gamma$, and
then show that the matching on $A$ restricts to a matching on $A\setminus B$.

\begin{rem}\label{R:M2}
 The chains $\gamma\in\Gamma$ which occur in the Pieri formula for cohomology
 are those whose Pieri chains have every transposition marked and whose Monk
 chains (if $\gamma\in\Gamma'_{k,p-1}$, since $\Gamma'_{k,p}$ does not occur)
 consists of a single transposition. 
 Restricting the argument given here to such chains furnishes a correction to
 the proof of the Pieri-type formula in~\cite{manfsp}. 
\end{rem}

We introduce some notation.
Express chains $\gamma$ in $\Gamma$ as a list of transpositions labeling
covers in $\gamma$, underline the marked transpositions and separate, if
necessary, the Pieri- and Monk- chains with $|$.
Thus for $v=421536$, 
 \begin{equation}\label{Eq:Gp_52}
  \Bigr( \underline{(4,6)}, \underline{(1,6)}, (2,5)\mid (3,5),(5,6)\Bigl)\ \in\ 
  \Gamma'_{5,2}
 \end{equation}
is the concatenation of Pieri chain from $421536$ to $531624$ with a Monk chain
to $532641$, which we display
\[
 \Blue{421}\Magenta{5}\,\Blue{3}\Magenta{6}
 \xrightarrow{\,\underline{(4,6)}\,}
 \Magenta{4}\Blue{216\,3}\Magenta{5}
 \xrightarrow{\,\underline{(1,6)}\,}
 \Blue{5}\Magenta{2}\Blue{16}\,\Magenta{3}\Blue{4}
 \xrightarrow{\,(2,5)\,}
 \Blue{531624}
   \mid 
 \Blue{53}\Magenta{1}\Blue{6}\Magenta{2}\,\Blue{4}
 \xrightarrow{\,(3,5)\,}
 \Blue{5326}\Magenta{1\,4}
 \xrightarrow{\,(5,6)\,}
 \Blue{532641}\,.
\]
The second transposition must be marked by Condition (P4).
The initial transposition and any subsequent transpositions relevant to
condition (P4) constitute the {\em initial subchain} of a given chain. 

\subsection{The first matching}\label{S:Matching1}
Define subsets $A_P$ and $A_M$ of $\Gamma$.
\begin{enumerate}
 \item[$A_P$:] The chains $\gamma\in\Gamma$ whose Pieri chain ends in an unmarked
        transposition $(i,k)$ which can be moved into the Monk chain to obtain a
        valid chain $\mu(\gamma)$ in $\Gamma$.
        It is possible that $\gamma\in\Gamma_{k,p}$ or $\Gamma_{k,p-1}$, so that
        it has no Monk chain.
 \item[$A_M$:] The chains $\gamma\in\Gamma$ whose Monk chain begins with a
        transposition $(i,k)$ which can be moved into the Pieri chain and 
        {\it left unmarked} to obtain a valid chain $\pi(\gamma)$ in $\Gamma$.
\end{enumerate}

The paired chains $\gamma,\mu(\gamma)$ and $\gamma,\pi(\gamma)$ contribute
opposite sign to the sum~\eqref{E:esum}.

Observe that the chain~\eqref{Eq:Gp_52} does not lie in $A_P$: if we
try to move the transposition $(2,5)$ from its Pieri chain to is Monk
chain, then that chain becomes $(2,5)(3,5)(5,6)$, which violates the
condition~\eqref{Eq:Monk_Condition} for a Monk chain as $a_1=2<3=a_2$.
However, it does lie in $A_M$: since $(2,5)\prec(3,5)$ the
transposition $(3,5)$ may be moved into the Pieri chain and left
unmarked to produce a valid Pieri chain.
As claimed, these two chains contribute opposite signs in the
sum~\eqref{E:esum}; the first has sign $-1$ and the second has sign $(-1)^2=1$.

In the definition of $A_M$, the condition that $\pi(\gamma)$ lies in
$\Gamma$ excludes chains $\gamma\in\Gamma$ having one of two exceptional forms.
 \begin{enumerate}
  \item[(E1):\ ] $\gamma\in\Gamma'_{k,p-1}$ has Monk chain consisting of a single
    transposition $(j,k)$ which can be moved into its Pieri chain and left
    unmarked to create a valid chain $\pi(\gamma)\in \Gamma_{k-1,p-1}$, and thus
    $\pi(\gamma)$ is {\it not} a chain in $\Gamma$. 

  \item[(E2):\ ] $\gamma\in\Gamma'_{k,p}\cup\Gamma'_{k,p-1}$ has Monk chain
    beginning with $(j,k)$ and Pieri chain of the form 
  \[
     \underline{(i_1,k)},\underline{(i_2,k)},\dotsc,\underline{(i_r,k)}\qquad
     \mbox{with}\quad i_1>i_2>\dotsb>i_r>j\,.
  \]
    Here, $r=p$ if $\gamma\in\Gamma'_{k,p}$ and $r=p{-}1$ if
    $\gamma\in\Gamma'_{k,p-1}$. 
    If we tried to move the transposition $(j,k)$ into the Pieri chain, it must
    be marked, by Condition (P4).
 \end{enumerate}
These cases will be treated in the next section.
We show that $\pi$ and $\mu$ are inverses and they define a matching on
$A:=A_P\cup A_M$. 

\begin{lem}\label{L:bijection}
 The sets $A_P$ and $A_M$ are disjoint, and $\pi,\mu$ are bijections between them.
\end{lem}

\begin{proof}
 Note that $\mu(A_P)\subset A_M$ and $\pi(A_M)\subset A_P$, and $\pi,\mu$ are
 inverses.
 We need only show that $A_P\cap A_M=\emptyset$.

 Suppose that the Pieri chain of $\gamma\in\Gamma$ ends in an unmarked
 transposition $(i,k)$ and its Monk chain begins with a transposition $(j,k)$. 
 Note that $(j,k)$ can be moved into the Pieri chain of $\gamma$ and left
 unmarked to create a valid chain if and only if $i<j$, by Condition (P3).
 This chain will not lie in $\Gamma$ if $\gamma$ has the 
 exceptional form (E1).
 Similarly, $(i,k)$ can be moved into the Monk chain to create a valid chain if
 and only if $i>j$, by Theorem~\ref{xkmonk}.
 Thus $\gamma$ cannot simultaneously lie in both $A_P$ and $A_M$.
\end{proof}

\begin{lem}\label{L:characterize}
  The set $A$ consists of chains $\gamma\in\Gamma$ that do not have one of the 
  exceptional forms, and  
  either their Pieri chain ends in an unmarked transposition $(i,k)$, or
  else their Monk chain begins with a transposition $(j,k)$ (or both).
\end{lem}

\begin{proof}
 Suppose that $\gamma\in\Gamma$ does not have one of the exceptional forms
 and its Monk chain begins with $(j,k)$.
 Then this can be moved into the Pieri chain to create a valid chain in $\Gamma$
 unless that Pieri chain ends in an unmarked $(i,k)$ with $i>j$, but then
 $(i,k)$ can be moved into the Monk chain.
 If the Monk chain does not begin with a transposition of the form $(j,k)$ and
 its Pieri chain ends in an unmarked $(i,k)$, then this can be moved into its
 Monk chain to create a valid chain in $\Gamma$.
\end{proof}

\subsection{The second matching}\label{S:Matching2}
 This is done in four steps with the first and most involved step matching every
 chain in $\Gamma_{k,p}$ with a chain in one of $\Gamma_{k,p-1}$,
 $\Gamma'_{k,p}$, or $\Gamma'_{k,p-1}$.
 The next three steps pair some of the remaining chains.
 Note that $p\geq 2$, as the base of our induction is the Monk formula with
 $p=1$.\medskip

{\bf Step 1.}
 Let $\gamma \in \Gamma_{k,p}$.
 Recall that a Pieri chain in $\Gamma_{k,p}$ is a chain in the $k$-Bruhat order
 with $p$ marked covers satisfying conditions (P1), (P2), (P3), and (P4) of
 Definition~\ref{D:Pieri_Chain}. 
 If no transposition $(k,\cdot)$ appears in $\gamma$, then $\gamma$ is also a
 chain in $\Gamma_{k-1,p}$, and we pair these two copies of $\gamma$ which
 contribute opposite signs to the sum~\eqref{E:esum}.
 Every chain in $\Gamma_{k-1,p}$ that is lacking a transposition of the form
 $(\cdot,k)$ is paired in this step.

 Now suppose that $\gamma\in\Gamma_{k,p}$ has a transposition of 
 the form $(k,\cdot)$, and let $(k,l)$ be the first such transposition in
 $\gamma$.
 The chain $\gamma$ has the form
 \begin{equation}\label{E:form}
   \gamma\ =\ (\dotsc,(k,l),\{(i_t,l)\}_{t=1}^r,
              \{\dotsc,(k,m_t)\}_{t=1}^s,\dotsc)\,,
 \end{equation}
 where $r,s\geq 0$; we also assume that all transpositions $(k,\cdot)$ and all
 transpositions $(\cdot,l)$ after $(k,l)$ were displayed in (\ref{E:form}).
 The parts in braces need not occur in $\gamma$, and we have
 suppressed the markings.

\begin{ex}\label{Ex:52173846}
 Suppose that $k=5$, $p=4$, and $v=52173846$.
 Then
 \begin{multline*}
   \Magenta{5}\Blue{2173\,84}\Magenta{6}\ 
   \xrightarrow{\,\underline{(1,8)}\,}\ 
   \Blue{6217}\Magenta{3}\,\Blue{8}\Magenta{4}\Blue{5}\ 
   \xrightarrow{\,\underline{(5,7)}\,}\ 
   \Blue{6}\Magenta{2}\Blue{174\,8}\Magenta{3}\Blue{5}\ 
   \xrightarrow{\,(2,7)\,}\ 
   \Blue{63}\Magenta{1}\Blue{74\,8}\Magenta{2}\Blue{5}\ 
   \xrightarrow{\,\underline{(3,7)}\,}\ \\
   \Blue{632}\Magenta{7}\Blue{4}\,\Magenta{8}\Blue{15}\ 
   \xrightarrow{\,\underline{(4,6)}\,}\ 
   \Magenta{6}\Blue{2384}\,\Magenta{7}\Blue{15}\ 
   \xrightarrow{\,(1,6)\,}\ 
   \Blue{7238}\Magenta{4\,6}\Blue{15}\ 
   \xrightarrow{\,(5,6)\,}\  \Blue{73286415}\,
\end{multline*}
is a Pieri chain in $\Gamma_{5,4}$.
We first use the intertwining relations of Lemma~\ref{L:length2}(2)(a) to
replace $\underline{(5,7)}(2,7)\underline{(3,7)}$ with 
$(2,5)\underline{(3,5)}(5,7)$ (unmarking the $(5,7)$), and obtain the chain
\begin{multline*}
   \Magenta{5}\Blue{2173\,84}\Magenta{6}\ 
   \xrightarrow{\,\underline{(1,8)}\,}\ 
   \Blue{6}\Magenta{2}\Blue{17}\,\Magenta{3}\Blue{845}\ 
   \xrightarrow{\,(2,5)\,}\ 
   \Blue{63}\Magenta{1}\Blue{7}\,\Magenta{2}\Blue{845}\ 
   \xrightarrow{\,\underline{(3,5)}\,}\ 
   \Blue{6327}\Magenta{1}\,\Blue{8}\Magenta{4}\Blue{5}\
   \xrightarrow{\,(5,7)\,}\  \\
   \Blue{632}\Magenta{7}\Blue{4}\,\Magenta{8}\Blue{15}\ 
   \xrightarrow{\,\underline{(4,6)}\,}\ 
   \Magenta{6}\Blue{2384}\,\Magenta{7}\Blue{15}\ 
   \xrightarrow{\,(1,6)\,}\ 
   \Blue{7238}\Magenta{4\,6}\Blue{15}\ 
   \xrightarrow{\,(5,6)\,}\  \Blue{73286415}\,.
\end{multline*}
The transpositions $(2,5)$, $\underline{(3,5)}$, and $(5,7)$ may now be commuted
past $(4,6)$ and $(1,6)$. 
Placing  the unmarked $(5,7)$ into the Monk chain 
gives the following chain in $\Gamma'_{5,3}$:
\begin{multline*}
  \Magenta{5}\Blue{217\,384}\Magenta{6}\ 
   \xrightarrow{\,\underline{(1,8)}\,}\ 
   \Blue{621}\Magenta{7}\,\Blue{3}\Magenta{8}\Blue{45}\ 
   \xrightarrow{\,\underline{(4,6)}\,}\ 
   \Magenta{6}\Blue{218\,3}\Magenta{7}\Blue{45}\ 
   \xrightarrow{\,(1,6)\,}\ 
   \Blue{7}\Magenta{2}\Blue{18}\,\Magenta{3}\Blue{645}\
   \xrightarrow{\,(2,5)\,}\\ 
   \Blue{73}\Magenta{1}\Blue{8}\,\Magenta{2}\Blue{645}\ 
   \xrightarrow{\,\underline{(3,5)}\,}\
   \Blue{73281645}\ \raisebox{-3pt}{\rule{0.4pt}{15pt}}\ 
   \Blue{7328}\Magenta{1}\,\Blue{6}\Magenta{4}\Blue{5}\ 
   \xrightarrow{\,(5,7)\,}\ 
   \Blue{7238}\Magenta{4\,6}\Blue{15}\ 
   \xrightarrow{\,(5,6)\,}\  \Blue{73286415}\,.
\end{multline*}
\end{ex}

 Given a chain $\gamma$ of the form~\eqref{E:form}, we transform it as
 in example~\ref{Ex:52173846} by first using the intertwining relations of
 Lemma~\ref{L:length2}(2)(a) to 
 \begin{equation}\label{E:intertw}
   \textrm{replace}\ \ (k,l),(i_1,l),\dotsc,(i_r,l)\ \ 
   \textrm{with}\ \ (i_1,k),\dotsc,(i_r,k),(k,l)\,.
 \end{equation}
 Let $\gamma''$ be the transformed chain. We then move  
  all transpositions now involving $k$ to the end of the chain
 using the commutation relations of Lemma~\ref{L:length2}(1).
 This is indeed possible as the following conditionds hold.

\begin{enumerate}
 \item[(i)] There are no transpositions $(\cdot,k)$ in $\gamma''$ other than those
       indicated in~\eqref{E:intertw}, as $\gamma$ is a chain in the $k$-Bruhat
       order. 
 \item[(ii)] Each transposition $(k,m_t)$ is the last transposition in $\gamma$
       involving $m_t$.
       Otherwise Condition (P3) would force $(k,m_t)$ to be
       marked, which is impossible by  Condition (P2), as $(k,l)$ precedes
       $(k,m_t)$ in $\gamma$. 
       Thus, any transpositions to the right of $(k,m_t)$ not involving $k$ 
       will commute with $(k,m_t)$.
 \item[(iii)] There are no transpositions $(i_t,\cdot)$ to the right of $(i_t,l)$ in
       $\gamma$. 
       Indeed, let $(i_t,l')$ be one.
       Then $l'<l$, by  Condition (P1). 
       We thus have the following subchain in $\gamma$:
 \[
     (k,l),\,\dotsc,\,(i_t,l),\,\dotsc,\,(i_t,l')\ .
 \]
      As $i_t<k<l'<l$, this subchain is forbidden by
      Lemma~\ref{L:Prohibition_1}.
\end{enumerate}

Let us split the chain obtained above from $\gamma''$ by commutations just before
the transposition $(k,l)$ to obtain the chain
\[
  \gamma'\ =\ (\dotsc,\,\{(i_t,k)\}_{t=1}^r\; \mid\;
    (k,l),\,\{(k,m_t)\}_{t=1}^s\;)\,.
\]
This is the concatenation of a chain in the ($k{-}1$)-Bruhat
order satisfying  Condition (P1) and a Monk chain. 
We describe how to mark $\gamma'$ to obtain a chain in either
$\Gamma'_{k,p}$ or $\Gamma'_{k,p-1}$.

If $\gamma$ does not begin with $(k,l)$ or if $\gamma$ begins with $(k,l)$ and
its initial subchain is just $(k,l)$, so that 
 \begin{equation}\label{spform}
   \gamma\ =\ ((k,l),\{(i_t,l)\}_{t=1}^r,
              \{(k,m_t)\}_{t=1}^s)\,,
 \end{equation}
then we mark the transpositions in $\gamma'$ that were marked in $\gamma$ and 
mark the transposition $(i_t,k)$ in $\gamma'$ if $(i_t,l)$ was marked in
$\gamma$. 
We also remove the mark (if any) from $(k,l)$.
This gives a valid marking of $\gamma'$.
Indeed, assume that $(i_t,l)$ was marked in $\gamma$.
By Condition (P2), no transposition $(i_t,\cdot)$ preceeds $(i_t,l)$ in
$\gamma$, and by (iii) above, no transposition $(i_t,\cdot)$ follows $(i_t,l)$,
so the transposition $(i_t,k)$ is the only transposition involving $i_t$ in
$\gamma'$ and so it may be marked and satisfy (P2).  
Since no transposition $(k,m_t)$ was marked in $\gamma$, we obtain a chain
$\gamma'$ in  $\Gamma'_{k,p-1}$ if $(k,l)$ was marked, and one in
$\Gamma'_{k,p}$ if $(k,l)$ was not marked. 
This chain $\gamma'$ contributes a sign opposite to that of 
$\gamma$ in the sum~\eqref{E:esum}.

If $\gamma$ begins with $(k,l)$ and its initial subchain contains other
transpositions, then the initial subchain of $\gamma'$ will involve
transpositions  $(j,m)$ with $m<l$ that are to the right of the initial  
subchain of $\gamma$.
To obtain a valid marking of $\gamma'$, first swap (in $\gamma$) the markings of
the last transpositions in the initial subchains of $\gamma$ and
$\gamma'$ (wherever the latter may appear in $\gamma$), and then proceed as above.
The chain $\gamma'$ will have $p{-}1$ marked
transpositions---losing the mark on $(k,l)$---except when the
initial subchain of $\gamma$ is just the transposition $(k,l)$ and the last
transposition in the initial subchain of $\gamma'$ is unmarked in $\gamma$. 
In that case, $\gamma'$ will have $p$ markings. 

We identify the image of
$\Gamma_{k,p}$ in each of $\Gamma_{k-1,p}$, $\Gamma'_{k,p}$, and
$\Gamma'_{k,p-1}$.
Call these images $\Gamma_{k-1,p}(1)$, $\Gamma'_{k,p}(1)$, and
$\Gamma'_{k,p-1}(1)$---the (1) indicates that these are the chains
paired in step 1 of this second matching.

 $\bullet$ $\Gamma_{k-1,p}(1)$: This is the intersection
 $\Gamma_{k,p}\cap\Gamma_{k-1,p}$ and it consists of those chains in 
 $\Gamma_{k-1,p}$ that do not contain a transposition of the form 
 $(\cdot,k)$.

 $\bullet$ $\Gamma'_{k,p}(1)$: 
 A chain $\gamma'$ is in $\Gamma'_{k,p}(1)$ if it is the image of a chain
 $\gamma\in\Gamma_{k,p}$ either whose first transposition $(k,l)$ involving $k$
 is unmarked, or else the initial subchain of $\gamma$ is just $(k,l)$ and the
 last transposition in the initial subchain of $\gamma'$ is unmarked in
 $\gamma$.
 In either case, $r=0$.
 The chain $\gamma'$ is obtained from $\gamma$ by trivially commuting all
 transpositions involving $k$ to the end with no intertwining.
 Thus, these are the chains in $\Gamma'_{k,p}$ whose Monk chain begins with
 $(k,l)$ and which have no transposition involving $k$ in their Pieri chain.

 $\bullet$ $\Gamma'_{k,p-1}(1)$:  This consists of the images of
 chains $\gamma$ in $\Gamma_{k,p}$ whose first transposition $(k,l)$
 involving $k$ is marked and, if $(k,l)$ is their initial subchain, then
 the last transposition in the initial subchain of $\gamma'$ is marked in
 $\gamma$.
 This is the set of chains in $\Gamma'_{k-1,p}$ such that the following hold.
\begin{enumerate}
 \item[(i)] The Monk chain begins with $(k,l)$.
 \item[(ii)] All transpositions
   $(i,k)$ in the Pieri chain intertwine with the transposition $(k,l)$ in (i)
   as in Lemma~\ref{L:length2}(2)(a). 
 \item[(iii)] If $(i,k)$ is the rightmost such transposition in the Pieri chain
   (necessarily the last transposition), then any other transposition $(i,m)$
   in the Pieri chain has $m>l$.
\end{enumerate}

The weakness of condition (iii) is explained by Lemma~\ref{L:technical},
which implies that if $(j,k)$ and $(j,m)$ are transpositions in the Pieri chain
of a chain in $\Gamma'_{k,p-1}(1)$, then $m>l$.

We describe the inverse transformations.
If $\gamma'\in\Gamma_{k-1,p}(1)$, then the inverse transformation
simply regards $\gamma'$ as a chain in $\Gamma_{k,p}$. 
If $\gamma'\in\Gamma'_{k,p}(1)$, then its Monk chain begins
with $(k,l)$ and its Pieri chain has no transpositions involving $k$.
Commute all transpositions $(k,\cdot)$ in its Monk chain back into its Pieri
chain as far left as possible to satisfy Condition (P1), preserving all
markings.
This gives a valid chain in $\Gamma_{k,p}$, except when $(k,l)$
becomes the initial subchain.
In that case, mark $(k,l)$ to satisfy Condition (P4) and unmark the
last transposition in the initial subchain of $\gamma'$ to obtain a chain in  
$\Gamma_{k,p}$.

Suppose that $\gamma'\in \Gamma'_{k,p-1}(1)$
satisfies conditions (i)--(iii) above.
By Lemma~\ref{L:technical}, if $(i,m)$ and $(i,k)$ are two 
transpositions in $\gamma'$, then $m>l$.
Thus, if we intertwine $(k,l)$ with all transpositions $(i_t,k)$ in
$\gamma'$ (the reverse of~\eqref{E:intertw}), we may commute the
transpositions $(k,l)$ and $(i_t,l)$ leftwards, as well as all
remaining transpositions $(k,m_t)$ in the Monk chain to obtain a 
chain $\gamma$ satisfying Condition (P1).
We then mark $(k,l)$.
This gives a valid chain in $\Gamma_{k,p}$ except possibly if
$(k,l)$ is its initial transposition without being its initial subchain, and
(\ref{spform}) holds. 
In that case, the last step in the transformation is to simply swap the markings
of the last transpositions of 
the initial subchains of $\gamma$ and $\gamma'$.

Lastly, we remark that the paired chains contribute opposite signs
to the sum~\eqref{E:esum}.\medskip

\noindent{\bf Step 2.}
 Let $\Gamma'_{k,p-1}(2.1)$ be the set of chains $\gamma$ in $\Gamma'_{k,p-1}$
 such that the following hold. 
\begin{enumerate}
 \item[(i)]   The Monk chain of $\gamma$ begins with a transposition $(k,l)$.
 \item[(ii)]  There is a transposition
              $(i,k)$ in the Pieri chain of $\gamma$ intertwining with the
              transposition $(k,l)$ as in 
              Lemma~\ref{L:length2}(2)(b), so that $(i,k)(k,l)$ becomes
              $(i,l)(i,k)$.  
 \item[(iii)] Let $(j,k)$ be the last transposition in the Pieri chain
              involving $k$.  
              By (ii) there is at least one.
              Then any other occurrence $(j,m)$ of $j$ has $m>l$.
\end{enumerate}
 These are chains in $\Gamma'_{k,p-1}$ that fail to be in
 $\Gamma'_{k,p-1}(1)$ only because of the way $(k,l)$ intertwines in (ii).

\begin{ex}
 Let $k=4$ and $p=3$.
 Then the following chain lies in $\Gamma'_{k,p-1}(2.1)$:
\begin{multline*}
  \Blue{3}\Magenta{1}\Blue{6\,5874}\Magenta{2}\ 
  \xrightarrow{\,\underline{(2,8)}\,}\ 
  \Blue{32}\Magenta{6}\,\Blue{5}\Magenta{8}\Blue{741}\ 
  \xrightarrow{\,(3,5)\,}\ 
  \Magenta{3}\Blue{28\,}\Magenta{5}\Blue{6741}\ 
  \xrightarrow{\,\underline{(1,4)}\,}\ 
  \Blue{5}\Magenta{2}\Blue{8\,}\Magenta{3}\Blue{6741}\ 
  \xrightarrow{\,(2,4)\,}
  \Blue{53826741}\ \raisebox{-3pt}{\rule{0.4pt}{15pt}}\\
  \Blue{538}\Magenta{2}\,\Blue{67}\Magenta{4}\Blue{1}\ 
  \xrightarrow{\,(4,7)\,}\ 
  \Blue{538}\Magenta{4\,6}\Blue{721}\ 
  \xrightarrow{\,(4,5)\,}\ 
  \Blue{53864721}\ .
\end{multline*}
Here, $l=7$, as the Monk chain begins with the transposition $(4,7)$.
If we try to intertwine $(4,7)$ with the transpositions in the Pieri
chain involving $k=4$, we obtain
\begin{eqnarray*}
  &&\Magenta{3}\Blue{28}\Magenta{5}\Blue{6741}\ 
   \xrightarrow{\,(1,4)\,}\ 
   \Blue{5}\Magenta{2}\Blue{8}\Magenta{3}\Blue{6741}\ 
   \xrightarrow{\,(2,4)\,}
   \Blue{538}\Magenta{2}\Blue{67}\Magenta{4}\Blue{1}\ 
   \xrightarrow{\,(4,7)\,}\ 
   \Blue{53846721}   \\ 
  &&\Magenta{3}\Blue{28}\Magenta{5}\Blue{6741}\ 
   \xrightarrow{\,(1,4)\,}\ 
   \Blue{528}\Magenta{3}\Blue{67}\Magenta{4}\Blue{1}\ 
   \xrightarrow{\,(4,7)\,}
   \Blue{5}\Magenta{2}\Blue{8467}\Magenta{3}\Blue{1}\ 
   \xrightarrow{\,(2,7)\,}\ 
   \Blue{53846721}   \\
  &&\Magenta{3}\Blue{28567}\Magenta{4}\Blue{1}\ 
   \xrightarrow{\,(1,7)\,}\ 
   \Magenta{4}\Blue{28}\Magenta{5}\Blue{6731}\ 
   \xrightarrow{\,(1,4)\,}
   \Blue{5}\Magenta{2}\Blue{8467}\Magenta{3}\Blue{1}\ 
   \xrightarrow{\,(2,7)\,}\ 
   \Blue{53846721}\ .
\end{eqnarray*}
 So while $(2,4)(4,7)$ intertwines as in Lemma~\ref{L:length2}(2)(a) to
 become $(4,7)(2,7)$, the next pair $(1,4)(4,7)$  intertwines as in
 Lemma~\ref{L:length2}(2)(b) to become $(1,7)(1,4)$.

 After this intertwining, we commute the transposition $(1,4)$ with
 $(2,7)$, then commute $(3,5)$ past $(1,7)(2,7)$, and place $(1,4)$
 in the Monk chain to obtain
\begin{multline}\label{Eq:Gamma(2.2)}
  \Blue{3}\Magenta{1}\Blue{6\,5874}\Magenta{2}\ 
  \xrightarrow{\,\underline{(2,8)}\,}\ 
  \Magenta{3}\Blue{26\,587}\Magenta{4}\Blue{1}
  \xrightarrow{\,\underline{(1,7)}\,}\ 
  \Blue{4}\Magenta{2}\Blue{6\,587}\Magenta{3}\Blue{1}
  \xrightarrow{\,(2,7)\,}\ 
  \Blue{43}\Magenta{6}\,\Blue{5}\Magenta{8}\Blue{721}\ 
  \xrightarrow{\,(3,5)\,}
  \Blue{43856721}\ \raisebox{-3pt}{\rule{0.4pt}{15pt}}\\ 
  \Magenta{4}\Blue{38}\Magenta{5}\,\Blue{6721}
  \xrightarrow{\,(1,4)\,}\ 
  \Blue{538}\Magenta{4\,6}\Blue{721}\ 
  \xrightarrow{\,(4,5)\,}\ 
  \Blue{53864721}\ .
\end{multline}
 We let $(1,7)$ inherit the mark on the original $(1,4)$.
 The result is a chain in $\Gamma'_{k,p-1}$, but one whose Monk chain
 begins with a transposition of the form $(i,k)=(1,4)$.
\end{ex}

 A chain $\gamma$ in $\Gamma'_{k,p-1}(2.1)$ has the following form
\[
  \gamma\ =\ (\ldots,\,\{(i_t,k)\}_{t=0}^r\ \mid\ 
              (k,l),\,\{(k,m_t)\}_{t=1}^s)\,,
\]
 where $r,s\geq 0$ and $(i_0,k)$ is the transposition $(i,k)$ of 
 condition (ii) above.
 Note that $l>m_1>\dotsb>m_t$ by Theorem~\ref{xkmonk}.
 We produce a new chain $\gamma'$ in $\Gamma'_{k,p-1}$ by first intertwining
 $(k,l)$ with the displayed transpositions $(i_t,k)$ to its left, and then
 commuting the transposition $(i_0,k)$ with the transpositions $(i_t,l)$ to its
 right as follows: 
 \begin{equation}\label{Eq:commuting}
  \begin{array}{rcl}
    \Blue{(i_0,k)},\,\dotsc,\,(i_r,k)\mid(k,l) 
        &=&\Blue{(i_0,k)},\,(k,l),\,(i_1,l),\,\dotsc,\,(i_r,l)\\
        &=&(i_0,l),\,\Blue{(i_0,k)},\,(i_1,l),\,\dotsc,\,(i_r,l)\rule{0pt}{15pt}\\
        &=&(i_0,l),\,(i_1,l),\,\dotsc,\,(i_r,l)\mid\Blue{(i_0,k)}\,.\rule{0pt}{15pt}
 \end{array}
\end{equation}
 By Lemma~\ref{L:technical}, each transposition $(i_t,l)$ may be 
 commuted leftwards in the Pieri chain of $\gamma$ to obtain a new chain
 $\gamma'$ satisfying Condition (P1). 
 As indicated, we let the Monk chain of $\gamma'$  begin with $(i_0,k)$ and
 declare the rest to be the Pieri chain.
 
\begin{rem}\label{R:M3}
 The chains $\gamma\in\Gamma'_{k,p-1}(2.1)$ which occur in the Pieri formula for
 cohomology have a Monk chain consisting only of $(k,l)$.
 In the proof given on page 94 of~\cite{manfsp}, it was assumed that $(i_r,k)$
 intertwines with $(k,l)$ as in Lemma~\ref{L:length2}(2)(b).
 As we see here, it may be the case that some other $(i_0,k)$ intertwines with
 $(k,l)$ in this manner.
\end{rem}

 We now mark the transpositions in $\gamma'$ that were marked in $\gamma$, and
 let $(i_t,l)$ inherit the mark of 
 $(i_t,k)$. This gives $p{-}1$ marks.
 If $(i_0,l)$ is the initial transposition of $\gamma'$ and the Pieri chain of
 $\gamma$ is not $((i_0,k),\ldots,(i_r,k))$, we need an extra step 
 to ensure that Condition (P4) holds. 
 First swap (in $\gamma$) the markings of the last transposition in the initial
 subchain of $\gamma$ and the transposition in $\gamma$ that will become the
 last transposition in the initial subchain of $\gamma'$, and then proceed as
 above.   
 We claim that this gives $\gamma'$ a valid marking, and therefore produces a
 chain in $\Gamma'_{k,p-1}$.

 Indeed, the only way that this could fail to be valid would be if
 the rightmost transposition $(j,l)$ in $\gamma$ involving $l$ was unmarked and
 had $j>i_0$, for then $(j,l)$ and $(i_0,l)$ would be adjacent in $\gamma'$ and
 Condition (P3) would force $(j,l)$ to be marked in $\gamma'$.
 But this gives a subchain of $\gamma'$ 
\[
    (j,l),\,(i_0,l),\,\dotsc,\,(i_0,k) \qquad\mbox{with}\quad i_0\,<\,j\,,
\]
 which is forbidden by Lemma~\ref{L:Prohibition_1}, since it is a chain in the
 ($k{-}1$)-Bruhat order, and it satisfies Conditions (P0) and (P1). Indeed, if
 $(h,k)$ is to the left of $(i_0,k)$ in $\gamma'$, then it was to the left of
 $(i_0,k)$ in $\gamma$; this means that $h$ cannot appear before $(h,k)$ in
 $\gamma'$ if $h>i_0$. 

 Because the transformation $\gamma \to \gamma'$ involves converting a
 transposition $(k,l)$ in the Monk chain of $\gamma$ into a transposition
 $(i_0,k)$ in the Monk chain of $\gamma'$, the two chains contribute opposite
 signs to the the sum~\eqref{E:esum}.

 Let $\Gamma'_{k,p-1}(2.2)$ be the set of chains $\gamma'$ obtained in this way
 from chains $\gamma$ in $\Gamma'_{k,p-1}(2.1)$ 
 This is the set of chains in $\Gamma'_{k,p-1}$ such that the following hold.
\begin{enumerate}
 \item[(i)] The Monk chain has a unique transposition of the form $(i,k)$.

 \item[(ii)] If the Pieri chain ends in an unmarked transposition $(j,k)$,
       then $j<i$. 

 \item[(iii)] The Pieri chain contains a transposition of the form
       $(i,\cdot)$. For any such transposition $(i,l)$, and for any
       transposition $(k,m)$ in the Monk chain, we have $l>m$. 

\end{enumerate}
 Note that the chain~\eqref{Eq:Gamma(2.2)} lies in
 $\Gamma'_{k,p-1}(2.2)$, for $k=4$ and $p=3$.
 The reason for condition (ii) is that if $(j,k)$ is an unmarked
 transposition just to the left of $(i_0,k)$ in $\gamma$, then $j<i_0$.
 Given a chain $\gamma'\in\Gamma'_{k,p-1}(2.2)$, we reverse the above procedure 
 to produce a chain in $\Gamma'_{k,p-1}(2.1)$.

 To see that it is indeed possible to reverse the above procedure, let
 $\gamma'\in\Gamma'_{k,p-1}(2.2)$. 
 Let $\delta'$ be its Pieri chain.
 By Conditions (i) and (ii), we may move the transposition $(i,k)$ from the Monk
 chain of $\gamma$ into its Pieri chain $\delta'$ to obtain a chain
 $\delta:=(\delta',(i,k))$ in the $(k{-}1)$-Bruhat order which satisfies
 Conditions (P0) and (P1).
 Let $(i,l)$ be the last transposition $(i,\cdot)$ in $\delta'$, and look at the
 transpositions beginning with $(i,l)$ in $\delta$.
 (Write $i_0=i$.)
\[
  (i_0,l),\, (i_1,l),\,\dotsc,\,(i_r,l),\,\dotsc,\,(i_0,k)\,.
\]
 The first step in this reverse procedure is to commute the transpositions
 $(i_j,l)$ to the end of the chain.
 We may do this, as the only obstruction would be if there were a transposition
 $(i_j,m)$ to the right of $(i_j,l)$ in $\delta$, and this is impossible.
 By condition (P1), we have $k\leq m<l$, and the subchain 
 $(i_0,l),\dotsc,(i_j,l),\dotsc,(i_j,m),\dotsc,(i_0,k)$ is forbidden by
 Lemma~\ref{L:quad}.

 Since the commutation is possible, we then just reverse the procedure
 in~\eqref{Eq:commuting} and reverse the marking procedure to obtain a valid
 chain in  $\Gamma'_{k,p-1}(2.1)$.\medskip

\noindent{\bf Step 3.}
 In this step, we pair some chains in $\Gamma'_{k,p-1}$ with chains in
 $\Gamma_{k-1,p}$ and $\Gamma'_{k,p}$.
 Define $\Gamma'_{k,p-1}(3)$ to be those chains in $\Gamma'_{k,p-1}$ such that the following hold.
\begin{enumerate}
 \item[(i)] The Monk chain has a unique transposition of the form $(i,k)$.

 \item[(ii)] If the Pieri chain ends in an unmarked transposition $(j,k)$,
       then $j<i$. 

 \item[(iii)] The Pieri chain has no transposition involving $i$.

\end{enumerate}

\begin{ex}
 For example, when $p=3$ and $k=4$, here is a chain in
 $\Gamma'_{k,p-1}(3)$.
\begin{multline*}
 \Blue{1}\Magenta{4}\Blue{2\,36}\Magenta{5}\ 
 \xrightarrow{\,\underline{(2,6)}\,}\ 
 \Blue{1}\Magenta{5}\Blue{2\,3}\Magenta{6}\Blue{4}\ 
 \xrightarrow{\,(2,5)\,}\ 
 \Blue{16}\Magenta{2\,3}\Blue{54}\ 
 \xrightarrow{\,\underline{(3,4)}\,}\ 
 \Blue{163254}\ \raisebox{-3pt}{\rule{0.4pt}{15pt}}\\ 
 \Magenta{1}\Blue{63}\Magenta{2}\,\Blue{54}\ 
 \xrightarrow{\,(1,4)\,}\ 
 \Blue{263}\Magenta{1\,5}\Blue{4}\ 
 \xrightarrow{\,(4,5)\,}\ 
 \Blue{263514}\,.
\end{multline*}
 Observe that we may move the initial transposition of its Monk chain
 into its Pieri chain, and mark it to obtain a valid chain in
 $\Gamma'_{k,p}$. 
\begin{multline*}
 \Blue{1}\Magenta{4}\Blue{2\,36}\Magenta{5}\ 
 \xrightarrow{\,\underline{(2,6)}\,}\ 
 \Blue{1}\Magenta{5}\Blue{2\,3}\Magenta{6}\Blue{4}\ 
 \xrightarrow{\,(2,5)\,}\ 
 \Blue{16}\Magenta{2\,3}\Blue{54}\ 
 \xrightarrow{\,\underline{(3,4)}\,}\ 
 \Magenta{1}\Blue{63\,}\Magenta{2}\Blue{54}\ 
 \xrightarrow{\,\underline{(1,4)}\,}\ 
 \Blue{263154}\ \raisebox{-3pt}{\rule{0.4pt}{15pt}}\\  
 \Blue{263}\Magenta{1\,5}\Blue{4}\ 
 \xrightarrow{\,(4,5)\,}\ 
 \Blue{263514}\,.
\end{multline*}
\end{ex}

 Given a chain $\gamma\in\Gamma'_{k,p-1}(3)$, we produce a chain
 $\gamma'$ by moving the transposition $(i,k)$ into its Pieri chain
 and then marking it. 
 This does not violate Conditions (P2) and (P3), by (ii) and (iii) above.
 If the Monk chain consists solely of $(i,k)$, then we obtain a chain in
 $\Gamma_{k-1,p}$, and otherwise a chain in $\Gamma'_{k,p}$.
 These images are characterized below.

 $\bullet$ $\Gamma_{k-1,p}(3)$ consists of chains in $\Gamma_{k-1,p}$ that end in a marked
 $(i,k)$; in other words, their inverse images are chains of the
 exceptional form (E1) in Section \ref{S:Matching1}.  

 $\bullet$ $\Gamma'_{k,p}(3)$ consists of those chains in $\Gamma'_{k,p}$ whose Pieri
 chain ends in a marked $(i,k)$ and whose Monk chain begins with a
 transposition $(k,l)$.

 The reverse procedure moves the marked $(i,k)$ into the Monk chain, and the
 paired chains contribute different signs to the sum~\eqref{E:esum}.\medskip

\noindent{\bf Step 4.}
 This involves the remaining chains having exceptional form (E2) in Section \ref{S:Matching1}.
 Let $\Gamma'_{k,p}(4)$ be those chains in $\Gamma'_{k,p}$ having exceptional
 form (E2) and $\Gamma'_{k,p-1}(4)$  be those chains in $\Gamma'_{k,p-1}$ having
 exceptional form (E2) and more than one transposition of the form $(j,k)$ in
 their Monk chain.
 (The chains in $\Gamma'_{k,p-1}$ having exceptional form (E2) and a single 
 transposition of the form $(j,k)$ in their Monk chain lie in
 $\Gamma'_{k,p-1}(3)$.)
 The matching between chains in $\Gamma'_{k,p}(4)$ on the left and 
 chains in  $\Gamma'_{k,p-1}(4)$ on the right is given below:
\[
  \Bigl(\underline{(i_1,k)},\,\dotsc,\,\underline{(i_p,k)}\,\mid\,
     (j,k)\,\dotsc\Bigr)\ \ \longleftrightarrow\ \
    \Bigl(\underline{(i_1,k)},\,\dotsc,\,\underline{(i_{p-1},k)}\,\mid\,
     (i_p,k),\,(j,k)\,\dotsc\Bigr)\,.
\]
 Here $i_1>i_2>\dotsb>i_p>j$ and, by (P4), all transpositions in both Pieri
 chains are marked. 
 Here is an example of such a chain in $\Gamma'_{6,2}(4)$:
\begin{multline*}
 \qquad
 \Blue{2361}\Magenta{4\,5}\Blue{7}\ 
 \xrightarrow{\,\underline{(5,6)}\,}\ 
 \Blue{2}\Magenta{3}\Blue{615}\,\Magenta{4}\Blue{7}\ 
 \xrightarrow{\,\underline{(2,6)}\,}\ 
 \Blue{2461537}\ \raisebox{-3pt}{\rule{0.4pt}{15pt}}\  \\
 \Magenta{2}\Blue{4615}\Magenta{3}\,\Blue{7}\ 
 \xrightarrow{\,(1,6)\,}\ 
 \Blue{346}\Magenta{1}\Blue{5}\Magenta{2}\,\Blue{7}\ 
 \xrightarrow{\,(4,6)\,}\ 
 \Blue{34625}\Magenta{1\,7}\ 
 \xrightarrow{\,(6,7)\,}\ 
 \Blue{3462571}\,.
 \qquad
\end{multline*}
 The matching moves the marked transposition $(1,6)$ into its Monk
 chain, giving a chain in  $\Gamma'_{6,3}(4)$
\begin{multline*}
 \qquad
 \Blue{2361}\Magenta{4\,5}\Blue{7}\ 
 \xrightarrow{\,\underline{(5,6)}\,}\ 
 \Blue{2}\Magenta{3}\Blue{615}\,\Magenta{4}\Blue{7}\ 
 \xrightarrow{\,\underline{(2,6)}\,}\ 
 \Magenta{2}\Blue{4615}\,\Magenta{3}\Blue{7}\ 
 \xrightarrow{\,\underline{(1,6)}\,}\ 
 \Blue{3461527}\ \raisebox{-3pt}{\rule{0.4pt}{15pt}}\ \\ 
 \Blue{346}\Magenta{1}\Blue{5}\Magenta{2}\,\Blue{7}\ 
 \xrightarrow{\,(4,6)\,}\ 
 \Blue{34625}\Magenta{1\,7}\ 
 \xrightarrow{\,(6,7)\,}\ 
 \Blue{3462571}\,.
 \qquad
\end{multline*}

 Note that the sets 
 $\Gamma_{k-1,p}(1)$, $\Gamma_{k-1,p}(3)$, $\Gamma'_{k,p}(1)$,
 $\Gamma'_{k,p}(3)$, $\Gamma'_{k,p}(4)$, $\Gamma'_{k,p-1}(1)$, 
 $\Gamma'_{k,p-1}(2.1)$, $\Gamma'_{k,p-1}(2.2)$, $\Gamma'_{k,p-1}(3)$, and
 $\Gamma'_{k,p-1}(4)$ are all disjoint. 
 Let $B$ be the union of these sets and $\Gamma_{k,p}$.

\subsection{Patching the matchings}\label{S:Matching3}

 We show that the two matchings (on the sets $A$ and $B$ defined in Sections \ref{S:Matching1} and \ref{S:Matching2}, respectively) include all chains in $\Gamma$, and that the
 matching on the set $A$ restricts to a matching on  $\Gamma\setminus B=A\setminus B$.
 Thus, we may patch the matching on $B$ with the matching
 on $\Gamma\setminus B$ to obtain a matching on $\Gamma$, which establishes the
 formula~\ref{E:esum}, and completes the proof of Theorem~\ref{kthm}.
 
\begin{lem}\label{L:union}
   $A\cup B = \Gamma$.
\end{lem}

\begin{proof}
 We have $\Gamma_{k,p}\subset B$ by definition.

 By Lemma~\ref{L:characterize}, $\Gamma_{k-1,p}(1)\cup\Gamma_{k-1,p}(3)$ 
 is the complement of $A$ in $\Gamma_{k-1,p}$.
 Similarly, $\Gamma'_{k,p}(1)\cup\Gamma'_{k,p}(3)\cup\Gamma'_{k,p}(4)$ is the
 complement of $A$ in $\Gamma'_{k,p}$.

 We consider $\Gamma'_{k,p-1}$.
 First note that the union $\Gamma'_{k,p-1}(1)\cup\Gamma'_{k,p-1}(2.1)$ consists
 of those chains $\gamma$ in $B$ whose Monk chain begins with $(k,l)$ and which
 furthermore satisfy the following condition.
 \begin{enumerate}
   \item[(iii$'$)] If the Pieri chain of $\gamma$ ends in $(i,k)$, then any other
              transposition $(i,m)$ in the Pieri chain has $m>l$.
 \end{enumerate}
 Thus if $\gamma\in\Gamma'_{k,p-1}$ has Monk chain beginning with $(k,l)$, 
 it lies in $A$ unless its Pieri chain does not end in an unmarked $(i,k)$.
 But this implies that it satisfies (iii$'$) above trivially.

 If the Monk chain of $\gamma\in\Gamma'_{k,p-1}\setminus B$ begins with $(j,k)$,
 then Lemma~\ref{L:characterize} implies that $\gamma\in A$, as the
 exceptional forms (E1) and (E2) of Section~\ref{S:Matching1} are chains in $B$.
\end{proof}

 \begin{lem}
   The matching on $A$ restricts to a matching $\Gamma\setminus B$.
 \end{lem}

\begin{proof}
 Since the matching on chains in $A$ does not change their number of marked
 covers, we consider this separately on  
 $\Gamma_{k,p}\cup\Gamma_{k-1,p}\cup\Gamma'_{k,p}$ and 
 $\Gamma'_{k,p-1}$.
 In the proof of Lemma~\ref{L:union} we showed that 
\[
  \bigl(\Gamma_{k,p}\cup \Gamma_{k-1,p}\cup\Gamma'_{k,p}\bigr)\setminus B
    \ =\ 
  \bigl(\Gamma_{k,p}\cup\Gamma_{k-1,p}\cup\Gamma'_{k,p}\bigr)\cap A\,.
\]
 This implies that the matching on $A$ restricts to a matching on this set.

 We show that the matching on $A$ restricts to a matching
 on $\Gamma'_{k,p-1}\cap A\cap B$, which implies that it restricts to a matching
 on $\bigl(\Gamma'_{k,p-1}\cap A\bigr)\setminus B$.
 First recall that $\Gamma'_{k,p-1}(1)\cup\Gamma'_{k,p-1}(2.1)$ is the set of all
 chains $\gamma\in\Gamma'_{k,p-1}$ whose Monk chain begins with $(k,l)$ and
 which satisfy Condition (iii$'$) in the proof of Lemma~\ref{L:union}.
 Also note that $\Gamma'_{k,p-1}(2.2)\cup\Gamma'_{k,p-1}(3)$ is the set of all
 chains $\gamma\in\Gamma'_{k,p-1}$ whose Monk chain has a unique transposition
 $(i,k)$ and which satisfy the following conditions. 
\begin{enumerate}
 \item[(ii)] If the Pieri chain ends in an unmarked transposition $(j,k)$,
       then $j<i$. 
 
 \item[(iii$''$)] If $(k,l)$ is in the Monk chain of $\gamma$ and $(i,m)$ in its
        Pieri chain, then $l<m$.
\end{enumerate}
 Note that $\Gamma'_{k,p-1}(4)\cap A=\emptyset$, as $A$ does not include chains
 having form (E2).

 Let $\gamma\in\Gamma'_{k,p-1}\cap A\cap B$.
 If the Monk chain of $\gamma$ begins with $(k,l)$, then its Pieri chain 
 ends in an unmarked $(i,k)$.
 This is moved into the Monk chain in $\mu(\gamma)$, which now has a unique
 transposition of the form $(\cdot,k)$.
 This new chain $\mu(\gamma)$ clearly satisfies (ii),
 and it satisfies (iii$''$), as $\gamma$ satisfies (iii$'$).
 It is not exceptional, as it has the form $\mu(\gamma)$.
 Thus $\mu(\gamma)\in\Gamma'_{k,p-1}\cap A\cap B$.

 On the other hand, if the Monk chain of $\gamma$ begins with $(i,k)$, then it
 has a single transposition of the form $(\cdot,k)$. 
 Condition (ii) implies that $\gamma\in A_M$, and $\pi(\gamma)$ moves the
 transposition $(i,k)$ into the Pieri chain.
 Then $\pi(\gamma)\in \Gamma'_{k,p-1}$, as $\gamma$ does not have one of the
 exceptional forms.
 Thus, the chain $\pi(\gamma)$ has a non-empty Monk chain that begins with
 $(k,l)$; furthermore, it satisfies (iii$'$), as $\gamma$ satisfies (iii$''$).
 This completes the proof. 
\end{proof}

\section{Related results}

There is a similar Pieri-type formula for the product of
a Grothendieck polynomial with $\fg_{(p)}(x_1,\ldots,x_k)$. 
This can be deduced from Theorem~\ref{kthm} by applying the
standard involution on the flag manifold $Fl_n$ which interchanges the Schubert
varieties $X_w$ and $X_{\omega_0w\omega_0}$. 
This involution induces an automorphism on $K^0(Fl_n)$ mapping the Schubert
class represented by $\fg_w(x)$ to the class represented by
$\fg_{\omega_0w\omega_0}(x)$. 
In particular, it maps $\fg_{(1^p)}(x_1,\ldots,x_k)$ to
$\fg_{(p)}(x_1,\ldots,x_{n-k})$. 
This involution maps the $k$-Bruhat order to the $(n-k)$-Bruhat order, and 
the order $\prec$ on labels of covers to the order
$\triangleleft$ defined by  
 \begin{equation}
  (a,b)\triangleleft (c,d)\qquad\mbox{if and only if}\quad
  (a<c)\spa\mbox{or}\spa (a=c\spandsp b> d)\,. 
 \end{equation}

\begin{thm}\label{kthm1} 
  We have that
 \begin{equation}
   \fg_v(x)\,\fg_{(p)}(x_1,\ldots,x_k)\ =\ 
   \sum_{(\gamma,\alpha)}(-1)^{\ell(\gamma)-p}\,\fg_{\eg}(x)\,,
 \end{equation}
where the sum is over all saturated chains in the $k$-Bruhat order (on $S_\infty$) 
\[
    \gamma\ \colon\ v=v_0\ \xrightarrow{\,({a_1,b_1})\,}\ 
                      v_1\ \xrightarrow{\,({a_2,b_2})\,}\ 
                      \cdots\ 
                   \xrightarrow{\,({a_q,b_q})\,}\ =\ \eg\,,\qquad 
          q=\ell(\gamma)\,,
\]
together with $p$ marked covers satisfying
\begin{enumerate}
 \item[(P1$'$)] $a_1\le a_2\le\dotsb\le a_q\,$.
 \item[(P2$'$)] if $(a_i,b_i)$ is marked, then $b_j\ne b_i$ for $j<i\,$.
 \item[(P3$'$)] if $(a_i,b_i)$ is unmarked and $i+1\le q$, then 
                 $(a_i,b_i)\triangleleft(a_{i+1},b_{i+1})$.
 \item[(P4$'$)] if $a_1=\dotsb=a_r$ and $b_1<\dotsb<b_r$ for some 
    $r\ge 1$, then $(a_r,b_r)$ is marked.
\end{enumerate}
 This formula has no cancellations.
\end{thm}

The special (and trivial) case of this when $k=1$ was given in Corollary~5.2
of~\cite{LRS04}.
There are versions of Corollary \ref{compressedformula} and
Theorem~\ref{T:unique} corresponding to the multiplication by 
$\fg_{(p)}(x_1,\ldots,x_k)$. 
These follow from the original ones above, so we omit them.

The Pieri-type formula of Theorem~\ref{kthm} is a common generalization of the
Pieri formula for Schubert polynomials in Theorem~ \ref{schub} and the Monk-type
formula in Theorem~\ref{kmonk}. 
The latter case is the specialization $p=1$, when the corresponding Pieri
chains have their first cover marked and all the other covers unmarked. 
On the other hand, the Monk-type formula can be
rearranged so that it is based on the order $\triangleleft$ on transpositions,
rather than $\prec$.
Then it becomes a special case of the Pieri-type formula in
Theorem~\ref{kthm1}.

Different special cases of the Pieri-type formulas above are the Pieri-type
formulas for Grothendieck polynomials corresponding to Grassmannian
permutations, which were obtained in~\cite{lencak}. 
We define some notation to state these formulas.
Given a partition $\lambda=(\lambda_1\ge \lambda_2\ge\dotsb\ge\lambda_k\ge 0)$,
call $|\lambda|:=\lambda_1+\dotsb+\lambda_k$ its {\em weight}.  
Let $r(\mu/\lambda)$ and $c(\mu/\lambda)$ denote the numbers of nonempty rows
and columns of a skew Young diagram $\mu/\lambda$. 
A skew diagram is a {\em horizontal} (respectively {\em vertical}) 
{\em strip} if it has no two boxes in the same column (respectively row).  

\begin{thm}\cite{lencak}\label{pieri1}
 Let $\lambda$ be a partition with at most $k$ parts.
 \begin{enumerate}
  \item [(1)] 
   ${\displaystyle \fg_\lambda(x_1,\ldots,x_k)\,\fg_{(p)}(x_1,\ldots,x_k)\ =\ 
   \sum_\mu
   (-1)^{|\mu|-|\lambda|}\,\binom{r(\mu/\lambda)-1}{|\mu/\lambda|-p}
    \,\fg_\mu(x_1,\ldots,x_k)}$,\newline
   where the sum is over all partitions $\mu$ with at most $k$ parts such that
   $\mu/\lambda$ is a horizontal strip of weight at least $p$.

  \item [(2)] 
    Suppose that $p<k$.
    Then
\[
     \fg_\lambda(x_1,\ldots,x_k)\,\fg_{(1^p)}(x_1,\ldots,x_k)\ =\ 
     \sum_\mu (-1)^{|\mu|-|\lambda|}\,
     \binom{c(\mu/\lambda)-1}{|\mu/\lambda|-p}\,\fg_\mu(x_1,\ldots,x_k)\,,
\]
    where the sum ranges over all partitions $\mu$ with at most $k$ parts such that
   $\mu/\lambda$ is a vertical strip of weight at least $p$.
 \end{enumerate}
\end{thm}

The first formula follows from the Pieri-type formula in Theorem \ref{kthm1},
and the second from the formula in Theorem \ref{kthm}. 
Indeed, given a Grassmannian permutation with descent in position $k$, the
corresponding chains in Theorem \ref{kthm1} are concatenations of subchains of
the following form (using the notation in Section 3), for different values of
$a$: 
\[ 
   ((a,b),\,(a,b{+}1),\,\ldots\,,(a,b{+}r))\,.
\]
Thus, by Condition (P3$'$), the transpositions $(a,b),\ldots,(a,b+r-1)$ must be
marked.
If this subchain is the initial one, then $(a,b+r)$ must also be marked, by
(P4$'$). 
Condition (P1$'$) guarantees that the entries $b+i$ corresponding to different
subchains are distinct, so Condition (P2$'$) is fulfilled. 
Applying the transpositions in such a chain to a
Grassmannian permutation corresponds to adding a horizontal strip to its
diagram, where each subchain contributes a row in the strip. 
We are free to choose the labels on the last transposition in each
subchain except the first---this explains the binomial coefficient in the first
formula. 
The second formula is similar.


\providecommand{\bysame}{\leavevmode\hbox to3em{\hrulefill}\thinspace}
\providecommand{\MR}{\relax\ifhmode\unskip\space\fi MR }
\providecommand{\MRhref}[2]{%
  \href{http://www.ams.org/mathscinet-getitem?mr=#1}{#2}
}
\providecommand{\href}[2]{#2}

\end{document}